\documentclass[10pt,twoside]{siamart1116}

\usepackage[english]{babel}
\usepackage{graphicx,epstopdf,epsfig}
\usepackage{amsfonts,epsfig,fancyhdr,graphics, hyperref,amsmath,amssymb}

\usepackage{mathdots}
\usepackage{algorithm}
\usepackage{algpseudocode}

\newcommand{\HE}{Name of Handling Editor}
\newcommand{\DoS}{Month/Day/Year}
\newcommand{\DoA}{Month/Day/Year}
\newcommand{\CA}{Heike Fa\ss bender}

\newcommand{\Names}{Peter Benner, Heike Fa\ss bender, Michel-Niklas Senn}
\newcommand{\Title}{The Hamiltonian Extended Krylov Subspace Method}

\newtheorem{remark}[theorem]{Remark}

\renewcommand{\theequation}{\thesection.\arabic{equation}}
\renewtheorem{theorem}{Theorem}[section]

\begin{document}

\bibliographystyle{plain}

\setcounter{page}{1}

\thispagestyle{empty}

 \title{\Title\thanks{Received
 by the editors on \DoS.
 Accepted for publication on \DoA.
 Handling Editor: \HE. Corresponding Author: \CA}}

\author{
Peter Benner\thanks{Max Planck Institute for Dynamics of Complex Technical Systems,
Sandtorstr. 1, 39106 Magdeburg, Germany (benner@mpi-magdeburg.mpg.de)}
\and
Heike Fa{\ss}bender\thanks{
Institute for Numerical Analysis, TU Braunschweig, Universit\"atsplatz 2, 38106 Braunschweig, Germany (h.fassbender@tu-braunschweig.de, m.senn@tu-braunschweig.de)}
\and
Michel-Niklas Senn\footnotemark[3]}

\markboth{\Names}{\Title}

\maketitle

\begin{abstract}
An algorithm for constructing a $J$-orthogonal basis of the extended Krylov subspace
$\mathcal{K}_{r,s}=\operatorname{range}\{u,Hu, H^2u,$ $ \ldots, $ $H^{2r-1}u, H^{-1}u, H^{-2}u, \ldots, H^{-2s}u\},$
where $H \in \mathbb{R}^{2n \times 2n}$ is a large (and sparse) Hamiltonian matrix is derived (for $r = s+1$ or $r=s$). Surprisingly, this allows for short recurrences involving at most five previously generated basis vectors. Projecting $H$ onto the subspace $\mathcal{K}_{r,s}$ yields a small Hamiltonian matrix. The resulting HEKS algorithm may be used in order to approximate $f(H)u$ where $f$ is a function which maps the Hamiltonian matrix $H$ to, e.g.,  a (skew-)Hamiltonian or symplectic matrix. Numerical experiments illustrate that approximating $f(H)u$ with the HEKS algorithm is competitive for some functions compared to the use of other (structure-preserving) Krylov subspace methods.
\end{abstract}



\section{Introduction}
Let $H\in \mathbb{R}^{2n \times 2n}$ be a nonsingular (large-scale) Hamiltonian matrix, that is  $J_nH=(J_nH)^T,$ where
$J_n = \left[\begin{smallmatrix} 0 & I_n \\ -I_n & 0\end{smallmatrix}\right]\in \mathbb{R}^{2n \times 2n}$
and $I_n$ is the $n \times n$ identity matrix.
We are interested in
computing a $J$-orthogonal basis of the extended Krylov subspace
\begin{align}\label{eq:extKrylov}
\mathcal{K}_{r,s} &:= \mathcal{K}_{2r}(H,u) + \mathcal{K}_{2s}(H^{-1},H^{-1}u)
= \operatorname{range}\{u,Hu, H^2u,\ldots, H^{2r-1}u, H^{-1}u, H^{-2}u, \ldots, H^{-2s}u\},
\end{align}
where  $u\in \mathbb{R}^{2n}$
and either $r = s +1$ or $r = s.$
That is, assuming
\[
\operatorname{dim } \mathcal{K}_{2r}(H,u)=2r \quad\text{ and }\quad \operatorname{dim } \mathcal{K}_{2s}(H^{-1},H^{-1}u) =2s,
\]
we are looking for a matrix $S_{r+s} \in \mathbb{R}^{2n \times 2(r+s)}$  with $J$-orthonormal columns ($S_{r+s}^TJ_nS_{r+s} = J_{r+s}$)  such that the columns of $S_{r+s}$ span the same subspace as
$\mathcal{K}_{2r}(H,u)+ \mathcal{K}_{2s}(H^{-1},H^{-1}u).$

Extended Krylov subspaces
\[
\operatorname{range}\{ b, A^{-1}b, Ab, A^{-2}b, A^2b, \ldots, A^{-k}b,A^kb\} = \mathcal{K}_k(A,b) + \mathcal{K}_k(A^{-1},A^{-1}b)
\]
for general nonsingular matrices $A \in \mathbb{C}^{n \times n}$ and a vector $b \in \mathbb{C}^n$ have been used for the numerical approximation of $f(A)b$ for a function $f$ and a large matrix $A$ at least since the late 1990s mainly inspired by \cite{DruK98,KniS10}. 
In case an orthogonal matrix $V$ has been constructed such that $\operatorname{range}(V)= \mathcal{K}_k(A,b) + \mathcal{K}_k(A^{-1},A^{-1}b),$ an approximation to $f(A)b$ can be obtained as
\begin{equation}\label{eq:proj}
f(A)b \approx Vf(V^TAV)V^Tb.
\end{equation}
 More on functions of matrices, the computation of $f(A)b$  and the approximation of $f(a)b$  via Krylov subspace methods can be found in the all-encompassing monograph \cite{Hig08}.

The idea of constructing a $J$-orthogonal basis for the extended Krylov subspace $\mathcal{K}_{r,s}$ \eqref{eq:extKrylov}
has first been considered in \cite{Mei11} in the context of approximating $\exp(H)u.$ The Hamiltonian Extended Krylov Subspace (short HEKS)  method 
presented in \cite{Mei11}  is a straightforward adaption of the algorithm for computing an orthogonal basis of an extended Krylov subspace described in \cite{KniS10}.
Our main finding in this paper is the observation that the HEKS algorithm allows for a short recurrence to generate $S_{r+s}.$

 We will explore the use of an $J$-orthogonal basis $S_{r+s}$ of the extended Krylov subspace $\mathcal{K}_{r,s}$ \eqref{eq:extKrylov} for approximating $f(H)u$ for a (large-scale) Hamiltonian matrix $H$ and a vector $u \in \mathbb{R}^{2n}.$
Following the idea from \eqref{eq:proj} we have
\[
f(H)u \approx  S_{r+s}f(H_{r+s}) J_{r+s}^TS_{r+s}^{T}J_nu
\]
where $H_{r+s} = J_{r+s}^TS_{r+s}^{T}J_nHS_{r+s}\in \mathbb{R}^{2(r+s) \times 2(r+s)}$ is a Hamiltonian matrix. That is, we can preserve the rich structural information inherent to the Hamiltonian structure of the matrix $H$. This would not be possible by computing a standard (orthogonal) basis $V \in \mathbb{R}^{2n \times 2(r+s)}$ of $\mathcal{K}_{r,s}$ as the matrix product $V^THV$ will in general not be a Hamiltonian matrix even if $H$ is Hamiltonian. Hence, the HEKS algorithm may be
used in particular in order to approximate $f(H)u$ where $f$ is a function which maps the Hamiltonian matrix $H$ to  a structured matrix such as  a (skew-)Hamiltonian or symplectic matrix.
Such a structure-preserving approximation of $f(H)u$ is, e.g., important in the context of symplectic exponential integrators for Hamiltonian systems, see, e.g., \cite{EirK19,HaiLW06,MeiW17,Mei11}.
A structure-preserving approximation of $f(H)u$ may also be computed using, e.g., an $J$-orthogonal basis $\tilde S_{2k}$ of the standard Krylov subspace $\operatorname{range} \{u, Hu, H^2u \ldots, H^{2k-1}u\}.$ Such a basis can be generated by the Hamiltonian Lanczos method \cite{BFS11,BenF97,W04}. Both approaches will be compared later on.

The paper is structured as follows: Section \ref{sec:prelim} summarizes some basic well-known facts about Hamiltonian and $J$-orthogonal matrices. In Section \ref{sec:idea} the general idea of generating the desired $J$-orthogonal basis $S_{r+s}$ of \eqref{eq:extKrylov} as proposed in \cite{Mei11}  is sketched. Then,  it is noted that the projected  matrices $H_{r+s} = J_{r+s}^TS_{r+s}^{T}J_nHS_{r+s} $ and $J_{r+s}^TS_{r+s}^{T}J_nH^{-1}S_{r+s} $ have  at most $10k,$ resp. $10k+2$, nonzero entries. The details are given in  Section \ref{sec:4} and in Section \ref{sec:5}.  The resulting efficient HEKS algorithm using short recursions is summarized in Section \ref{sec:alg}.
The rather long and technical constructive proof for our claim is deferred to the Appendix \ref{sec:derivation}.
In Section \ref{sec:numex} the approximation of $f(H)u$ using the HEKS algorithm is compared to the approximation by the extended Krylov subspace method \cite{KniS11} and by the Hamiltonian Lanczos method \cite{BenF97}.

\section{Preliminaries}\label{sec:prelim}
Here we list some properties of Hamiltonian and $J$-orthogonal matrices useful for the following discussion.

\begin{enumerate}
\item $J_n= \left[\begin{smallmatrix} 0 & I_n \\ -I_n & 0\end{smallmatrix}\right]\in \mathbb{R}^{2n \times 2n}$ is orthogonal and skew-symmetric, $J_n^T=J_n^{-1} =-J_n$.
 \item Let $H\in\mathbb{R}^{2n\times2n}.$
$H$ is Hamiltonian if and only if there exist matrices $E$, $B=B^T$, $C=C^T$ $\in\mathbb{R}^{n\times n}$ such that
    \begin{align*}
        H = \begin{bmatrix}
            E & B \\
            C & -E^T
        \end{bmatrix}.
    \end{align*}
    \item Let $H\in\mathbb{R}^{2n\times2n}$ be a nonsingular Hamiltonian matrix. Then $H^{-1}$ is Hamiltonian as well.
\item The eigenvalues of a Hamiltonian matrix $H$ occur
in pairs $\{\lambda, -\lambda\}$ if $\lambda$ is real or purely imaginary, or in quadruples $\{\lambda, \overline{\lambda}, -\lambda, -\overline{\lambda}\}$ otherwise. That is, the spectrum of
a Hamiltonian matrix is symmetric with respect to both the real and the imaginary
axis.
\item A matrix $S \in \mathbb{R}^{2n \times 2n}$ is called \emph{symplectic}  if  $S^TJ_nS=J_n$. Its columns are \emph{$J$-orthogonal}.
    \item Let $S\in\mathbb{R}^{2n\times2n}$ be a symplectic matrix. Then $S^{-1} = J_n^TS^TJ_n$ is symplectic as well.
    \item Let $H\in\mathbb{R}^{2n\times2n}$ be a Hamiltonian matrix and $S\in\mathbb{R}^{2n\times2n}$ be a  symplectic matrix. Then $S^{-1}HS$ is a Hamiltonian matrix.
\item
Let $S \in \mathbb{R}^{2n \times 2m}, m \leq n,$  have $J$-orthogonal columns, $S^TJ_nS=J_m.$
Let $H \in \mathbb{R}^{2n \times 2n}$ be Hamiltonian.
\begin{enumerate}
    \item The matrix $J_m^TS^TJ_n$ is the left inverse of $S$, $J_m^TS^TJ_nS = I_{2m}.$
    \item The matrix $(J_m^TS^TJ_n)HS$ is Hamiltonian.
\end{enumerate}
\end{enumerate}
Numerous further properties of the sets of these matrices (and their interplay) have been studied in the literature, see, e.g., \cite{MMT03} and the references therein. In particular, $J_n$ induces a skew-symmetric bilinear form $\langle\cdot,\cdot\rangle_{J_n}$ on $\mathbb{R}^{2n}$ defined by $\langle x,y \rangle_{J_n} = y^TJ_nx$ for $x,y\in \mathbb{R}^{2n}.$ Hamiltonian matrices are skew-adjoint with respect to the bilinear form $\langle\cdot,\cdot\rangle_{J_n}$, while symplectic matrices are orthogonal with respect to $\langle\cdot,\cdot\rangle_{J_n}$. The $2n \times 2n$ symplectic matrices form a Lie group, the $2n \times 2n$ Hamiltonian matrices the associated Lie algebra.

Assume that  a matrix  $S_k =[V_k~~W_k]\in \mathbb{R}^{2n\times 2k}$  with $J$-orthogonal columns is given with $V_k=[v_1~~v_2~~\cdots~~v_k]$ and $W_k=[w_1~~w_2~~\cdots~~w_k] \in \mathbb{R}^{2n \times k}$.
Two additional vectors $x, J_nx \in \mathbb{R}^{2n}$ can be added to $S_k$ to generate a matrix $S_{k+1}=[V_{k+1}~~W_{k+1}]\in \mathbb{R}^{2n\times 2k+2}$ with $J$-orthogonal columns by $J$-orthogonalizing the vectors $x$ and  $J_nx$ against all column vectors $v_j,w_j$ of $S_k$ via
\begin{align*}
    v_{k+1} &= x -S_{k}J_{k}^TS_{k}^TJ_nx,\\
    w_{k+1} &= (J_nv_{k+1}) -S_{k}J_{k}^TS_{k}^TJ_n(J_n v_{k+1}), \quad w_{k+1} = w_{k+1}/(v_{k+1}^TJ_nw_{k+1}).
\end{align*}

\section{Idea of the HEKS Algorithm}\label{sec:idea}

Let a Hamiltonian matrix $H \in \mathbb{R}^{2n \times 2n}$ and a vector $u_1\in \mathbb{R}^{2n},$ $\|u_1\|_2=2,$ be given.
 Assume that $\operatorname{dim } \mathcal{K}_{2r}(H,u_1)=2r$ and $\operatorname{dim } \mathcal{K}_{2s}(H^{-1},H^{-1}u_1) =2s$.
The goal is to construct a matrix $S_{r+s} \in \mathbb{R}^{2n \times 2(r+s)}$  with $J$-orthonormal columns ($S_{r+s}^TJ_nS_{r+s} = J_{r+s}$)  such that the columns of $S_{r+s}$ span the same subspace as
$\mathcal{K}_{2r}(H,u_1)+ \mathcal{K}_{2s}(H^{-1},H^{-1}u_1).$

In \cite{Mei11} it is suggested to construct the matrix $S_{r+s}$ in the following way (assuming that no breakdown occurs):
\begin{enumerate}
\item We start with the two vectors in $\mathcal{K}_{2}(H,u_1)$ and construct
\[
S_1 = \begin{bmatrix} u_1\mid v_1\end{bmatrix} \in \mathbb{R}^{2n \times 2}
\]
 with $S_1^TJ_nS_1 = J_1$  and $\operatorname{range}\{S_1\} = \mathcal{K}_{2}(H,u_1).$ This corresponds to the choice $r=1, s=0.$
\item Thereafter we take the two vectors in $\mathcal{K}_{2}(H^{-1},H^{-1}u_1)$ and construct
\[
S_2 = \begin{bmatrix}y_1&u_1 \mid x_1 &v_1\end{bmatrix} = \begin{bmatrix}Y_1 & U_1 \mid X_1 & V_1\end{bmatrix}\in \mathbb{R}^{2n \times 4}
\]
 with $S_2^TJ_nS_2 = J_2$  and $\operatorname{range}\{S_2\} = \mathcal{K}_{2}(H,u_1)+ \mathcal{K}_{2}(H^{-1},H^{-1}u_1).$
This corresponds to the choice $r=s=1.$
\end{enumerate}
We proceed in this fashion by alternating between the subspaces $\mathcal{K}_{2r}(H,u_1)$ and $\mathcal{K}_{2s}(H^{-1},H^{-1}u_1).$ Assume that a matrix
\[
S_{2k} =\begin{bmatrix}Y_{k} & U_{k} \mid X_{k} & V_{k}\end{bmatrix} \in \mathbb{R}^{2n \times 4k}, \qquad Y_k, U_k, X_k, V_k \in \mathbb{R}^{2n \times k}
\]
 with $J$-orthonormal columns has been constructed such that
its columns span the same space as $\mathcal{K}_{2k}(H,u_1)+ \mathcal{K}_{2k}(H^{-1},H^{-1}u_1).$
The following three steps are repeated until the desired symplectic basis has been generated:
\begin{itemize}
\item[(3)] Construct $u_{k+1}$ and $v_{k+1}$ and set
\begin{align*}
S_{2k+1} &= \begin{bmatrix}Y_{k} & U_{k} & u_{k+1} \mid X_{k} & V_{k} & v_{k+1}\end{bmatrix}
 = \begin{bmatrix}Y_{k} & U_{k+1} \mid X_{k} & V_{k+1}\end{bmatrix} \in \mathbb{R}^{2n \times 4k+2}
\end{align*}
 with
\[
 U_{k+1} =\begin{bmatrix} U_k & u_{k+1}\end{bmatrix}, V_{k+1}=\begin{bmatrix}V_k & v_{k+1}\end{bmatrix} \in \mathbb{R}^{2n \times k+1}
\]
 such that $S_{2k+1}^TJ_nS_{2k+1} = J_{2k+1}$  and $\operatorname{range}\{S_{2k+1}\} = \mathcal{K}_{2k+2}(H,u_1)+ \mathcal{K}_{2k}(H^{-1},H^{-1}u_1).$ The new vectors are added as the last column to the $U,$ resp. $V$-matrix.
\item[(4)] Construct $y_{k+1}$ and $x_{k+1}$ and set
\begin{align*}
S_{2k+2} &= \begin{bmatrix}y_{k+1} & Y_{k} & U_{k+1} \mid x_{k+1} & X_{k} & V_{k+1}\end{bmatrix}
= \begin{bmatrix}Y_{k+1} & U_{k+1} \mid X_{k+1} & V_{k+1}\end{bmatrix} \in \mathbb{R}^{2n \times 4k+4}
\end{align*}
 with
\[
Y_{k+1}=\begin{bmatrix}y_{k+1} & Y_k\end{bmatrix},  X_{k+1}=\begin{bmatrix}x_{k+1} & X_k\end{bmatrix} \in \mathbb{R}^{2n \times k+1}
\]
 such that $S_{2k+2}^TJ_nS_{2k+2} = J_{2k+2}$  and $\operatorname{range}\{S_{2k+2}\} = \mathcal{K}_{2k+2}(H,u_1)+ \mathcal{K}_{2k+2}(H^{-1},H^{-1}u_1).$ The new vectors are added as the first column to the $Y$, resp. $X$-matrix.
\item[(5)] Set $k = k+1.$
\end{itemize}
We refrain from restating the algorithm given in \cite{Mei11} which implements the approach stated above in a straightforward way using long recurrences.
As usual, a Krylov recurrence of the form
 \[HS_{2k} = S_{2k}H_{2k} + \text{ some rest term},  \]
for $r=s=k$ and
\[HS_{2k+1} = S_{2k+1}H_{2k+1} + \text{some rest term} \]
for $r=s+1=k$ holds,
where $H_{2k}= J_{2k}^TS_{2k}^{T}J_nHS_{2k} \in \mathbb{R}^{4k\times 4k}$ and $H_{2k+1} = J_{2k+1}^TS_{2k+1}^{T}J_nHS_{2k+1} \in \mathbb{R}^{4k+2 \times 4k+2}$ are Hamiltonian matrices.
In the next two sections we describe the very special forms of the projected  matrices $H_{2k}$ and
$H_{2k+1}$ as well as
$J_{2k}^TS_{2k}^{T}J_nH^{-1}S_{2k} $ and
$ J_{2k+1}^TS_{2k+1}^{T}J_nH^{-1}S_{2k+1}. $
These matrices have  at most $10k,$ resp. $10k+2$, nonzero entries. A constructive proof for our claim is given in Appendix \ref{sec:derivation}, while the resulting efficient HEKS algorithm using short recursions is summarized in Section \ref{sec:alg}.

\section{Projection $ J_{r+s}^TS_{r+s}^{T}J_nHS_{r+s} $ of the Hamiltonian matrix $H$}\label{sec:4}
Assume that
\[
S_{r+s} = \begin{bmatrix}Y_s & U_r \mid X_s & V_r\end{bmatrix}, \qquad Y_s, X_s \in \mathbb{R}^{2n \times s},
\quad U_r, V_r \in \mathbb{R}^{2n \times r},
\]
with $J$-orthogonal columns has been constructed with the HEKS algorithm (as before, we assume that $r=s$ or $r = s+1$).
Then the projected  Hamiltonian matrix
\[
H_{r+s} = J_{r+s}^TS_{r+s}^{T}J_nHS_{r+s} \in \mathbb{R}^{2(r+s) \times 2(r+s)}
\]
has a very special form with at most $2r+8s$ nonzero entries.
Let us first note that
\[
H_{r+s} = {\small
\begin{bmatrix}
-X_s^TJ_nHY_s & -X_s^TJ_nHU_r & -X_s^TJ_nHX_s & -X_s^TJ_nHV_r\\
-V_r^TJ_nHY_s & -V_r^TJ_nHU_r & -V_r^TJ_nHX_s & -V_r^TJ_nHV_r\\
Y_s^TJ_nHY_s & Y_s^TJ_nHU_r & Y_s^TJ_nHX_s & Y_s^TJ_nHV_r\\
U_r^TJ_nHY_s & U_r^TJ_nHU_r & U_r^TJ_nHX_s & U_r^TJ_nHV_r
\end{bmatrix}},
\]
where the blocks are of size either $s\times s,$ $r\times r,$ $s\times r,$ or $r\times s.$ As will be proven in Appendix \ref{sec:derivation},
ten of these blocks are zero, three are diagonal (denoted by $\Delta_{s}, \Theta_{r}, \Lambda_{s}$), one symmetric tridiagonal (denoted by $T_{r}$)
and two anti-bidiagonal (denoted by $B_{sr}$), i.e.,
\begin{equation}\label{eq_struct1}
H_{r+s} =
\begin{bmatrix}
0 & 0 & \Lambda_{s} & B_{sr}\\
0 & 0 & B_{sr}^T & T_{r}\\
\Delta_{s} & 0 & 0 & 0\\
0 & \Theta_{r} &0 & 0
\end{bmatrix}
\end{equation}
with
\begin{center}
\begin{minipage}[t]{0.3\textwidth}
\begin{align*}
\Delta _s &= \operatorname{diag}(\delta_s, \ldots, \delta_1) \in \mathbb{R}^{s \times s}, \\
\Theta_r &=\operatorname{diag}(\vartheta_1, \ldots, \vartheta_r) \in \mathbb{R}^{r \times r}, \\
\Lambda_s &=\operatorname{diag}(\lambda_s, \ldots, \lambda_1) \in \mathbb{R}^{s \times s}, \\
\end{align*}
\end{minipage}\qquad
\begin{minipage}[t]{0.35\textwidth}
\begin{align*}
T_r &= \begin{bmatrix}
 \alpha_1 &\beta_2&& \\
\beta_2 & \ddots &\ddots& \\
 & \ddots &\ddots& \beta_r\\
 &  &\beta_r& \alpha_r
\end{bmatrix}\in \mathbb{R}^{r \times r}, \\
\end{align*}
\end{minipage}
\end{center}
and either
\[
B_{r-1,r} =\begin{bmatrix}
&&&\gamma_{r-1} & \mu_{r}\\
 &  &\iddots &  \mu_{r-1}& \\
 & \iddots & \iddots &  & \\
\gamma_{1} & \mu_{2} &  & &
\end{bmatrix}\in \mathbb{R}^{r-1\times r} \quad \text{~~ if ~~} r =s+1,
\]
or
\[
B_{rr} =\begin{bmatrix}
&&&\gamma_{r} \\
 &  &\iddots &  \mu_{r} \\
 & \iddots & \iddots &   \\
\gamma_{1} & \mu_{2} &
\end{bmatrix}\in \mathbb{R}^{r\times r} \quad \text{~~ if ~~} r =s.
\]
In particular, it holds for $j = 1, \ldots, s$
\begin{align*}
\delta_j &= y_j^TJ_nHy_j,\qquad
\lambda_j = -x_j^TJ_nHx_j,
\end{align*}
and for $j = 1, \ldots, r$
\begin{align*}
\vartheta_j &= u_j^TJ_nHu_j,\qquad
\alpha_j = -v_j^TJ_nHv_j,\qquad
\gamma_j = -x_j^TJ_nHv_j,
\end{align*}
and for $j=2,\ldots, r$
\begin{align*}
\beta_j &= -v_j^TJ_nHv_{j-1}\qquad
&\mu_j &= -x_{j-1}^TJ_nHv_j.
\end{align*}

We summarize this in the following theorem.
\begin{theorem}\label{theo1}
Let $H \in \mathbb{R}^{2n \times 2n}$ be a Hamiltonian matrix. Let $r+s = n$ and either $r = s +1$ or $r=s.$ Then
in case the procedure described in Section \ref{sec:idea} does not break down for
$u_1 \in \mathbb{R}^{2n}$ with $\|u_1\|_2 = 1$ there exists a symplectic matrix $S\in \mathbb{R}^{2n \times 2n}$ such that $Se_{s+1}=u_1,$
\[
\operatorname{range}\{S\} = \mathcal{K}_{2r}(H,u_1)+ \mathcal{K}_{2s}(H^{-1},H^{-1}u_1),
\]
and
\[
S^{-1}HS = H_{r+s}
\]
with $H_{r+s}=H_n$ as in \eqref{eq_struct1}.
\end{theorem}
\begin{proof}
A constructive proof is given in Section \ref{sec:derivation}.
\end{proof}

\begin{remark}\label{rem}
In case the Hamiltonian matrix $H$ can be written in the form $H=JK$ with the symmetric matrix $K$ and $K$ is positive definit, all inner products of the form $w^TJHw$ and $w^TJH^{-1}w$ are negative,
as $w^TJHw = w^TJJKw = -w^TKw <0$ and as with $K$ its inverse is symmetric and positive definite. Thus, in this case, all
$\delta_j$ and $\vartheta_j$ are negative, while  all $\lambda_j$ and $\alpha_j$ are positive. Such Hamiltonian matrices have been considered in \cite{Amo03,Amo06}.
\end{remark}

Theorem \ref{theo1} implies
\begin{align*}
H \begin{bmatrix}Y_k& U_k&X_k &V_k \end{bmatrix} &=
S\left[\begin{array}{c|c||c|c}
0 & 0 & 0 & 0\\
0 & 0 & \Lambda_{k} & B_{kk}\\\hline
0 & 0 & B_{kk}^T & T_{k}\\
0 & 0 & \mu_{k+1}e_1^T & \beta_{k+1}e_k^T\\
0 & 0 &0 & 0\\
\hline\hline
0 & 0 & 0 & 0\\
\Delta_{k} & 0 & 0 & 0\\\hline
0 & \Theta_{k} &0 & 0\\
0 & 0 &0 & 0\\
\end{array}\right].
\end{align*}
From this, we obtain the  HEKS-recursion for $r=s=k$
\begin{align}\label{eq:receven}
HS_{2k}
 &=S_{2k}H_{2k} + \mu_{k+1}u_{k+1}e_{2k+1}^T + \beta_{k+1}u_{k+1}e_{4k}^T,
\end{align}
%
while for $r=s+1=k+1$ we have
\begin{align}\label{eq:recodd}
HS_{2k+1}
 &=S_{2k+1}H_{2k+1} +  (\gamma_{k+1}y_{k+1}+ \beta_{k+2}u_{k+2})e_{4k+2}^T .
\end{align}

\section{Projection $ J_{r+s}^TS_{r+s}^{T}J_nH^{-1}S_{r+s} $ of the Hamiltonian matrix $H^{-1}$}\label{sec:5}
%
Assume that Theorem \ref{theo1} holds.
As $H_{n} = S^{-1}HS \in \mathbb{R}^{2n \times 2n} $ is Hamiltonian, its inverse $H_n^{-1} = S^{-1}H^{-1}S$ is Hamiltonian as well. Not only $H_n$ has a nice sparse structure \eqref{eq_struct1}, but also its inverse. From that we can derive the special forms of $ J_{2k}^TS_{2k}^{T}J_nH^{-1}S_{2k} $ and $ J_{2k+1}^TS_{2k+1}^{T}J_nH^{-1}S_{2k+1}.$

Let $S=S_n = \begin{bmatrix}Y_s  &U_r \mid X_s& V_r\end{bmatrix} \in \mathbb{R}^{2n \times 2n},$ $Y_s, X_s \in \mathbb{R}^{2n \times s},$ $U_r,V_r \in \mathbb{R}^{2n \times r},$ where $r+s=n$ and either $r=s$ or $r = s+1.$
Due to $S_n^{-1} = J_n^TS_n^TJ_n,$ we have
\begin{align*}
H_n^{-1} &={\small
\begin{bmatrix}
-X_s^TJ_nH^{-1}Y_s & -X_s^TJ_nH^{-1}U_r & -X_s^T J_nH^{-1}X_s & -X_s^TJ_nH^{-1}V_r\\
-V_r^TJ_nH^{-1}Y_s & -V_r^TJ_nH^{-1}U_r & -V_r^T J_nH^{-1}X_s & -V_r^TJ_nH^{-1}V_r\\
Y_s^TJ_nH^{-1}Y_s & Y_s^TJ_nH^{-1}U_r & Y_s^T J_nH^{-1}X_s & Y_s^TJ_nH^{-1}V_r\\
U_r^TJ_nH^{-1}Y_s & U_r^TJ_nH^{-1}U_r & U_r^T J_nH^{-1}X_s & U_r^TJ_nH^{-1}V_r
\end{bmatrix}}\\
& = \begin{bmatrix}
0 & 0 & \Delta_{s}^{-1} & 0\\
0 & 0 & 0 & \Theta_{r}^{-1}\\
E_{s} & G_{sr} & 0 & 0\\
G_{sr}^T &  F_{r} & 0&0
\end{bmatrix}
\end{align*}
with $E_s \in \mathbb{R}^{s \times s},$ $F_r \in \mathbb{R}^{r \times r},$ $G_{sr}\in \mathbb{R}^{s \times r}$ such that
\[
\begin{bmatrix}
\Lambda_{s} & B_{sr}\\ B_{sr}^T & T_{r}
\end{bmatrix}
\begin{bmatrix}
E_{s} & G_{sr}\\
G_{sr}^T & F_{r}
\end{bmatrix} =I
\]
holds and $\Delta_{s}, \Theta_{r}, \Lambda_{s}, T_{r}, B_{sr}$ from \eqref{eq_struct1}.
The matrices $E_s,F_r$ and $G_{sr}$ have a special structure like $\Lambda_s, T_r$ and $B_{sr}$: $F_r$ is diagonal, $G_{sr}$ anti-bidiagonal as $B_{sr}$ and $E_{s}$ is symmetric tridiagonal;
\begin{align*}
F_r &=\operatorname{diag}(
f_{11}, f_{22}, \ldots, f_{rr}), 
\qquad
E_s =\begin{bmatrix}
e_{ss} & e_{s-1,s} &  && \\
e_{s-1,s} & \ddots & \ddots & & \\
 & \ddots&\ddots& e_{12} & \\
 &  &e_{12} & e_{11}
\end{bmatrix} = E_s^T , 
\end{align*}
and either
\[
G_{r-1,r} =\begin{bmatrix}
&&&g_{r-1,r-1} & g_{r-1,r}\\
 &  &\iddots &  g_{r-2,r-1}& \\
 & g_{22} & \iddots &  & \\
g_{11} & g_{12} &  & &
\end{bmatrix}\in \mathbb{R}^{r-1\times r} \quad \text{~~ if ~~} r =s+1,
\]
or
\[
G_{rr} =\begin{bmatrix}
&&&g_{rr} \\
 &  &\iddots &  g_{r-1,r} \\
 & g_{22} & \iddots &   \\
g_{11} & g_{12} &
\end{bmatrix}\in \mathbb{R}^{r\times r} \quad \text{~~ if ~~} r =s.
\]

Next, the projected matrices $J_{2k}^TS_{2k}^TJ_nH^{-1}S_{2k}$ and $J_{2k+1}^TS_{2k+1}^TJ_nH^{-1}S_{2k+1}$ will be described.
Let
\begin{align*}\label{eq:Tlj}
\mathfrak E_j = \begin{bmatrix}I_j\\0\end{bmatrix}\in \mathbb{R}^{r \times j},\qquad
 \mathfrak F_\ell = \begin{bmatrix}0\\I_\ell\end{bmatrix}\in \mathbb{R}^{s \times \ell},\qquad
\mathfrak T_{\ell j} = \begin{bmatrix}\mathfrak F_\ell\\ &\mathfrak E_j\\&&\mathfrak F_\ell \\&&&\mathfrak E_j
\end{bmatrix}\in \mathbb{R}^{2n \times 2(\ell+j)}
\end{align*}
for $j \leq r,$ $\ell \leq s.$ Thus, for $2k \leq n$ it holds
\begin{align*}\label{eq_SnTk}
S_n \mathfrak T_{kk} = S_{2k} \in \mathbb{R}^{2n \times 4k}\qquad \text{and} \qquad
S_n \mathfrak T_{k,k+1} = S_{2k+1} \in \mathbb{R}^{2n \times 4k+2},
\end{align*}
as well as
\begin{align*}
S_nJ_n  \mathfrak T_{kk} &= \begin{bmatrix}-X_k& -V_k&Y_k &U_k \end{bmatrix} =S_{2k}J_{2k} \in \mathbb{R}^{2n \times 4k},\\
S_nJ_n  \mathfrak T_{k,k+1} &= \begin{bmatrix}-X_k& -V_{k+1}&Y_k &U_{k+1} \end{bmatrix} =S_{2k+1} J_{2k+1}\in \mathbb{R}^{2n \times 4k+2}.
\end{align*}
Hence, we obtain
\begin{align}
J_{2k}^TS_{2k}^TJ_nH^{-1}S_{2k} &=
\mathfrak T_{kk}^T J_n^TS_n^TJ_nH^{-1}S_n\mathfrak T_{kk}  = \mathfrak T_{kk}^T H_n^{-1}\mathfrak T_{kk}  \nonumber
= \mathfrak T_{kk}^T   \begin{bmatrix}
0 & 0 & \Delta_{s}^{-1} & 0\\
0 & 0 & 0 & \Theta_{r}^{-1}\\
E_{s} & G_{sr} & 0 & 0\\
G_{sr}^T &  F_{r} & 0&0
\end{bmatrix}\mathfrak T_{kk} \nonumber\\
&=  \begin{bmatrix}
0 & 0 & \mathfrak F_k^T\Delta_{s}^{-1}\mathfrak F_k & 0\\
0 & 0 & 0 & \mathfrak E_k^T\Theta_{r}^{-1}\mathfrak E_k\\
\mathfrak F_k^TE_{s}\mathfrak F_k & \mathfrak F_k^TG_{sr}\mathfrak E_k & 0 & 0\\
\mathfrak E_k^TG_{sr}^T\mathfrak F_k &  \mathfrak E_k^TF_{r} \mathfrak E_k& 0&0
\end{bmatrix} 
=  \begin{bmatrix}
0 & 0 & \Delta_{k}^{-1} & 0\\
0 & 0 & 0 & \Theta_{k}^{-1}\\
E_{k} & G_{kk} & 0 & 0\\
G_{kk}^T&  F_{k}& 0&0
\end{bmatrix} \label{eq_myHinv}
\end{align}
and
\begin{align}
J_{2k+1}^TS_{2k+1}^TJ_nH^{-1}S_{2k+1} &=
\mathfrak T_{k,k+1}^T J_n^TS_n^TJ_nH^{-1}S_n\mathfrak T_{k,k+1}
=  \begin{bmatrix}
0 & 0 & \Delta_{k}^{-1} & 0\\
0 & 0 & 0 & \Theta_{k+1}^{-1}\\
E_{k} & G_{k,k+1} & 0 & 0\\
G_{k,k+1}^T&  F_{k+1}& 0&0
\end{bmatrix}. \label{eq_myHinv2}
\end{align}
The HEKS-recurrences for $H^{-1}$ are given by
\begin{align}\label{eq_receven}
H^{-1}S_{2k} = S_{2k}\left( J_{2k}^TS_{2k}^TJ_nH^{-1}S_{2k}\right) +(e_{2k,2k+1}x_{2k+1}+g_{2k,2k+1}v_{2k+1})e_1^T
\end{align}
and
\begin{align}\label{eq_recodd}
H^{-1}S_{2k+1} = S_{2k+1}\left( J_{2k+1}^TS_{2k+1}^TJ_nH^{-1}S_{2k+1}\right) +x_{2k+1}(e_{2k,2k+1}e_1^T+g_{2k+1,2k+1}e_{2k+1}^T).
\end{align}

\section{HEKS Algorithm}\label{sec:alg}
The HEKS algorithm is summarized in Fig. \ref{alg}. The algorithm as given assumes that no breakdown occurs. Clearly, any division by zero will result in a serious breakdown. As can be seen from \eqref{eq:receven} a lucky breakdown occurs in case
$\mu_{k+1}=\beta_{k+1}=0$  or $u_{k+1}=0$, as $\operatorname{range}\{S_{2k}\}=\mathcal{K}_{2k}(H,u_1)+ \mathcal{K}_{2k}(H^{-1},H^{-1}u_1)$ is $H$-invariant.
Moreover, \eqref{eq:recodd} shows that in case $\gamma_{k+1}=\beta_{k+2}=0,$ a lucky breakdown occurs, as $\operatorname{range}\{S_{2k+1}\}=\mathcal{K}_{2k+2}(H,u_1)+ \mathcal{K}_{2k}(H^{-1},H^{-1}u_1)$ is $H$-invariant. Similarly, lucky breakdown can be read off of \eqref{eq_receven} and \eqref{eq_recodd} resulting in an $H^{-1}$-invariant subspace.

In case the Hamiltonian matrix $H$
can be written in the form $H=JK$ with a symmetric positive definite matrix $K$, all inner products of the form $w^TJHw$ and $w^TJH^{-1}w$ are negative (see Remark \ref{rem}). Hence, most scalars by which is divided in Algorithm \ref{alg} are nonzero and do not cause breakdown.

 Implemented efficiently such that each matrix-vector product as well as each linear solve is computed only once, the algorithm requires (for adding $4$ vectors) in the for-loop
\begin{itemize}
\item $4$ matrix-vector-multiplications with $H$,
\item $3$ linear solves with $H$ (efficiently implemented in the form $(JH)x=(Jb)$ making use of the symmetry of $JH$),
\item $14$ scalar products.
\end{itemize}
 Any multiplication of a vector $w$ by $J_n$ should be implemented by rearranging the upper and the lower part of the vector $w$. That is, let $w =\left[\begin{smallmatrix}w_1\\w_2\end{smallmatrix}\right],$ then $J_n w = \left[\begin{smallmatrix}w_2\\-w_1\end{smallmatrix}\right].$

Without some form of re-$J$-orthogonalization the HEKS algorithm suffers from the same numerical difficulties as any other Krylov subspace method.


{\small
\begin{algorithm}[!ht]
    \caption{HEKS with short recurrences }
    \label{alg}
    \begin{algorithmic}[1]
        \Require Hamiltonian matrix $H{\in}\mathbb{R}^{2n\times2n}$, $u_1{\in}\mathbb{R}^{2n}$ with $\|u_1\|_2=1$
        \Ensure  a) $S_{2k}= [y_k~\cdots~y_1~u_1~\cdots~u_k\mid x_k~\cdots~x_1~v_1~\cdots~v_k\in \mathbb{R}^{2n \times 4k}$ with $S_{2k}^TJ_nS_{2k}=J_{2k}$ and $H_{2k} =  J_{2k}S_{2k}^TJ_nHS_{2k}$ as in \eqref{eq_struct1}\newline
b) parameters $\lambda_j, \delta_j, \alpha_j,\gamma_j, \vartheta_j$ for $j = 1, \ldots, k$ and $\beta_j, \mu_j$ for $j=2, \ldots, k$ which determine $H_{2k}$\newline
(for $S_{2k+1}\in \mathbb{R}^{2n \times 4k+2}$ the algorithm needs to be modified appropriately)

\State $u_1 = u_1/\|u_1\|_2$  \Comment{Set up $S_1 = [u_1\mid v_1]$}
\State $\vartheta_1 = u_1^TJ_nHu_1$
\State $v_1 = Hu_1/\vartheta_1$
\State $f_{11} = u_1^TJ_nH^{-1}u_1$\Comment{Set up $S_2 = [y_1~u_1\mid x_1~v_1]$}
\State $w_x = H^{-1}u_1 -f_{11}v_1$
\State $x_1 = w_x/\|w_x\|_2$
\State $y_1 = H^{-1}x_1/x_1^TJ_nH^{-1}x_1$
\State $\lambda_1 = -x_1^TJ_nHx_1$ and  $\delta_1 = y_1^TJ_nHy_1$
\State $\alpha_1 = -v_1^TJ_nHv_1$ and  $\gamma_1 = -x_1^TJ_nHv_1$\Comment{Set up $S_3 = [y_1~u_1~u_2\mid x_1~v_1~v_2]$}
\State $w_u = Hv_1-\gamma_1y_1-\alpha_1u_1$
\State $u_2 = w_u/\|w_u\|_2$
\State $\vartheta_2 = u_2^TJ_nHu_2$
\State $v_2 = Hu_2/\vartheta_2$
\State $e_{11} = y_1^TJ_nH^{-1}y_1$\Comment{Set up $S_4 = [y_2~y_1~u_1~u_2\mid x_2~x_1~v_1~v_2]$}
\State $ g_{11} = y_1^TJ_nH^{-1}u_1,$ and $g_{12}=y_1^TJ_nH^{-1}u_2$
\State $w_x = H^{-1}y_1 -e_{11}x_1 - g_{11}v_1-g_{12}v_2$
\State $x_2 = w_x/\|w_x\|_2$
\State $y_2 = H^{-1}x_2/(H^{-1}x_2)^TJ_nx_2$
\State $\lambda_2 = -x_2^TJ_nHx_2$ and $\delta_2 = y_2^TJ_nHy_2$
\For{$j=3, 4, \ldots, k$}
\State $\alpha_{j-1} = -v_{j-1}^TJ_nHv_{j-1}$ and  $\beta_{j-1} = -v_{j-1}^TJ_nHv_{j-2}$\Comment{Set up $S_{2j-1}$}
\State $\gamma_{j-1}= -x_{j-1}^TJ_nHv_{j-1}$ and $ \mu_{j-1}=-x_{j-2}^TJ_nHv_{j-1}$
\State $w_u = Hv_{j-1}-\gamma_{j-1} y_{j-1}-\mu_{j-1}y_{j-2}-\beta_{j-1} u_{j-2}-\alpha_{j-1} u_{j-1}$
\State $u_j = w_u/ \|w_u\|_2$
\State $\vartheta_j = u_j^TJ_nHu_j$
\State $v_j = Hu_j/\vartheta_j$
\State $g_{j-1,j-1}= y_{j-1}^TJ_nH^{-1}u_{j-1}$ and $g_{j-1,j}=y_{j-1}^TJ_nH^{-1}u_{j}$ \Comment{Set up $S_{2j}$}
\State $e_{j-1,j-1}= y_{j-1}^TJ_nH^{-1}y_{j-1}$ and $e_{j-2,j-1}=y_{j-1}^TJ_nH^{-1}y_{j-2}$
\State $w_x = H^{-1}y_{j-1} -e_{j-1,j-1}x_{j-1} -e_{j-2,j-1}x_{j-2}- g_{j-1,j-1}v_{j-1}-g_{j-1,j}v_j$
\State $x_j = w_x/\|w_x\|_2$
\State $y_j = H^{-1}x_j/(H^{-1}x_j)^TJ_nx_j$
\State $\lambda_j = -x_j^TJ_nHx_j$ and $\delta_j = y_j^TJ_nHy_j$
\EndFor
\State $\alpha_k = -v_k^TJ_nHv_k$ and $\beta_k = -v_k^TJ_nHv_{k-1}$
\State $\gamma_k =-x_k^TJ_nHv_k$ and $ \mu_k=-x_{k-1}^TJ_nHv_k$
    \end{algorithmic}
\end{algorithm}
}

\section{Numerical Experiments}\label{sec:numex}
In this section, we demonstrate experimentally that the HEKS algorithm may be useful   for approximating $f(H)u$ for a (large-scale) Hamiltonian matrix $H \in \mathbb{R}^{2n \times 2n}$ and a vector $u \in \mathbb{R}^{2n}, \|u\|_2 = 1,$ via
\begin{align}\label{eq_approxsymp}
f(H)u \approx \widetilde Sf(\widetilde H) J_{2\ell}^T\widetilde S^{T}J_nu
\end{align}
with the $2n \times 2\ell$ $J$-orthogonal matrix $\widetilde S$ and
 the $2\ell \times 2\ell $ Hamiltonian matrix  $\widetilde H = J_{2\ell}^T\widetilde S^{T}J_nH\widetilde S.$
We consider two methods to construct $\widetilde S$:
\begin{itemize}
\item the HEKS method (Algorithm \ref{alg}) which generates a $J$-orthogonal matrix $\widetilde S$ such that $\operatorname{range}(\widetilde S)=\mathcal{K}_{r,s}$ with $r=s=\frac{\ell}{2}$ or $r-1=s=\frac{\ell-1}{2},$ depending on whether $\ell$ is even or odd. Then $f(H)u$ can be  approximated via $\widetilde Sf(\widetilde H) e_{s+1}$ (as  due to the construction $J_{2\ell}^T\widetilde S^{T}J_nu=e_{s+1}$),
\item  the  Hamiltonian Lanczos method (HamL) \cite{BenF97,BFS11,W04} which generates a $J$-orthogonal matrix $\widetilde S$ such that
$\operatorname{range}(\widetilde S)=\mathcal{K}_{2\ell}(H,u).$ Then $f(H)u$ can be  approximated via $\widetilde Sf(\widetilde H) e_1$ (as  due to the construction $J_{2\ell}^T\widetilde S^{T}J_nu=e_1$).
\end{itemize}
These methods are compared to the corresponding unstructured methods
\begin{itemize}
\item the extended Krylov subspace method (EKSM) \cite{KniS10},
\item the standard Arnoldi method \cite{GolVL13},
\end{itemize}
which generate an orthogonal matrix $Q$ such that either $\operatorname{range}(Q)=  \mathcal{K}_{r,s}$ or $\operatorname{[range}(Q)=  \mathcal{K}_{2\ell}(H,u).$ Then $f(H)u$ can be  approximated via $Qf(Q^THQ)e_1$ (as by construction, $Q^Tu=e_1$ holds).

Only functions $f$ which map $H$ to a structured matrix are dealt with. In particular, we consider
\begin{itemize}
\item $f(H) = \exp(H):$  the exponential function of a Hamiltonian matrix is a symplectic matrix \cite{HarC06},
\item $f(H) = \cos(H):$  $\cos(H)$ is a skew-Hamiltonian matrix (as a sum of even powers of $H$),
\item $f(H) = \operatorname{sign}(H):$ $\operatorname{sign}(H)$ is a Hamiltonian matrix \cite{MMT06}.
The matrix sign function is defined for any matrix $X\in \mathbb{C}^{n \times n}$ having no pure imaginary eigenvalues by $\operatorname{sign}(X) = X(X^2)^{-\frac{1}{2}}$ \cite{Hig94,Hig08}. An equivalent definition is $\operatorname{sign}(X) = T\operatorname{diag}(-I_p,I_q)T^{-1},$ where the Jordan decomposition of $X = T\operatorname{diag}(J_1,J_2)T^{-1}$ is such that the $p$ eigenvalues of $J_1$ are assumed to lie in the open left half-plane, while the $q$ eigenvalues  of $J_2$ lie in the open right half-plane. The Newton iteration $S_0=X,$ $S_{k+1}=\frac{1}{2}(X_k+X_k^{-1})$ converges quadratically to $\operatorname{sign}(X)$ \cite{Rob71}.
\end{itemize}
Utilizing HEKS or HamL, the projected matrix $\widetilde H$ is Hamiltonian again, so that $f(\widetilde H)$ has the same structure as $f(H),$ while the projected matrix $Q^THQ$ as well as $f(Q^THQ)$ obtained via EKSM and Arnoldi have no particular structure.
Such a structure-preserving approximation of $f(H)u$ is, e.g., important in the context of symplectic exponential integrators for Hamiltonian systems, see, e.g., \cite{EirK19,HaiLW06,MeiW17,Mei11}.


All experiments are performed in MATLAB R2021b on an Intel(R) Core(TM) i7-8565U CPU @ 1.80GHz   1.99 GHz with 16GB RAM. Our MATLAB implementation employs the standard MATLAB function \texttt{expm} and \texttt{funm(H,@cos)} as well as \texttt{signm} from the Matrix Computation Toolbox \cite{Hig02}. The experimental code used to generate the results presented in the following subsection can be found at \cite{zenodo}. All algorithms are run to yield a $1000 \times 30$ matrix whose columns span the corresponding (extended) Krylov subspace. All methods are implemented  using full re-($J$)-orthogonalization.
The accuracy of the approximation for HEKS and HamL is measured in terms of the relative error $\|f(H)u-\widetilde Sf(\widetilde H) J_{2\ell}^T\widetilde S^{T}J_nu\|_2/\|f(H)u\|_2$, while $\|f(H)u-Qf(Q^THQ)Q^Tu\|_2/\|f(H)u\|_2$ is used for EKSM and Arnoldi.

\subsection{Example 1}
Inspired by \cite[Example 4.1]{KniS10}, our first test matrix is a diagonal Hamiltonian matrix $H_1 = \operatorname{diag} (D,-D)$ with a diagonal $500 \times 500$ real matrix $D$ whose eigenvalues are log-uniformly distributed in the interval $[10^{-1},1].$ EKSM will preserve the symmetry of $H$, while HEKS and HL will not.
\begin{figure}[ht]
 \begin{center}
\includegraphics[width=0.8\textwidth]{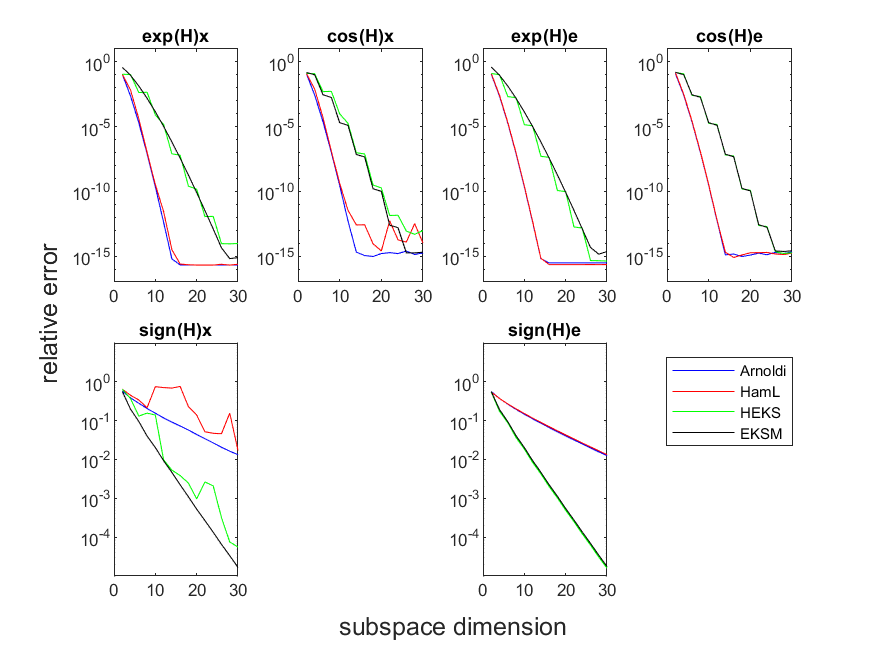}
 \caption{Diagonal Hamiltonian matrix $H_1 = \operatorname{diag}(A, -A)$ with \texttt{A = diag(logspace(-1,0,500));}, two different choices of the starting vector \texttt{x = randn(1000,1)} and \texttt{e = ones(1000,1)}}
 \end{center}
 \label{fig1}
 \end{figure}

In Fig. \ref{fig1}, the relative accuracy of all four methods is displayed, using a random starting vector $x$ (plots in the two leftmost columns) as well as a starting vector of all ones (plots in the two rightmost columns).
The Hamiltonian Lanczos method and the Arnoldi method perform alike just as the HEKS algorithm and the EKS method.
For the functions $\exp$ and $\cos$   the HEKS approximation makes significant progress only every other iteration step (that is, whenever the columns of $\widetilde S$ span $\mathcal{K}_{k,k-1}$). The same holds for the EKSM approximation of $\cos(H)x$ and $\cos(H)e$, but not for the approximation of $\exp(H)x$ and $\exp(H)e.$ The HEKS algorithm adds the vectors from $\mathcal{K}_{r,s}$ in a different order than EKSM: HEKS alternates between adding two vectors from $\mathcal{K}_{2k}(H,u)$ and adding two vectors from $\mathcal{K}_{2k}(H^{-1},H^{-1}u)$ , while EKSM alternates between adding one vector from $\mathcal{K}_{2k}(H,u)$ and adding one from $\mathcal{K}_{2k}(H^{-1},H^{-1}u)$ (for $u = x$ or $u=e$). Thus, the columns of $\widetilde S$ and $Q$ span the same subspace only every other step. Adding vectors from $\mathcal{K}_{2k}(H^{-1},H^{-1}u)$ does not seem to be relevant for the HEKS approximations $\exp(H)u$ and $\cos(H)u$ as well for the EKSM approximation of $\cos(H)u.$ For the EKSM approximation of $\exp(H)u$ some convergence progress can be observed in every iteration step, but the overall convergence is similar to that of the HEKS approximation. In summary, the use of an extended Krylov subspace does not improve the convergence for these examples compared to the approximations computed using the Arnoldi method or the Hamiltonian Lanczos method. The latter two methods converge about twice as fast as the first two.

But for the matrix sign function, the two methods based on the extended Krylov subspace converge faster than the other two. They do make progress in every iteration step. It is clearly beneficial to use an extended Krylov subspace here.

The HEKS algorithm requires  $34$ matrix-vector-multiplications with $H$,  $21$ linear solves with $H$ and $104$ scalar products to construct the $1000 \times 30$ matrix $\widetilde S.$ In contrast, the ESKM requires  $15$ matrix-vector-multiplications with $H$,  $14$ linear solves with $H$ and
 $493$ scalar products. As $H$ in this example is diagonal, the linear solves and matrix-vector multiplications require less arithmetic operations  than scalar products. Hence, the HEKS algorithm is faster than EKSM and requires less storage. Of course, the situation will change for more practically relevant
 examples with a more complex sparsity pattern. But it remains to
 note that there is a big difference in the number of scalar products
 to be performed, which is not due to the matrix structure but the
 difference of the short-term Lanczos-style and long-term Arnoldi-style
 recursions in the non-symmetric case.

\subsection{Example 2}
As a second example we use the Hamiltonian matrix
$ H_2 = \left[\begin{smallmatrix} A & -G\\ -Q & -A^T\end{smallmatrix}\right]\in \mathbb{R}^{1998 \times 1998}$
from Example 15 of  the collection of benchmark examples for the numerical solution of continuous-time algebraic Riccati equations \cite{BenLM95}.
The matrix has a complex spectrum with real and imaginary parts between $-2$ and $2.$

\begin{figure}[ht]
 \begin{center}
\includegraphics[width=0.8\textwidth]{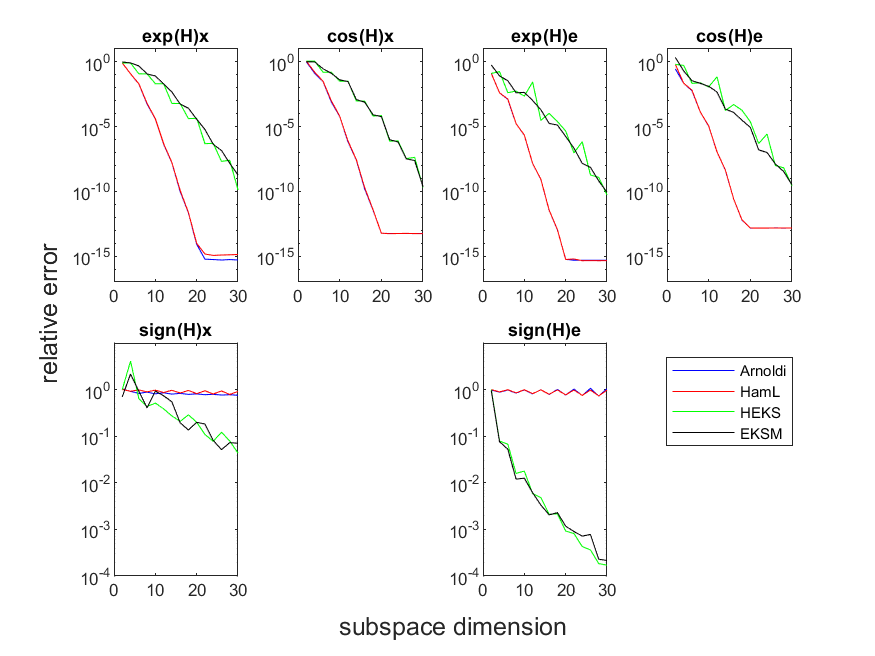}
 \caption{Hamiltonian matrix $H_2$, two different choices of the starting vector \texttt{x = randn(1998,1)} and \texttt{e = ones(1998,1)}}
 \end{center}
 \label{fig2}
 \end{figure}

Fig. \ref{fig2} provides the same information as in Fig. \ref{fig1}.  Our findings from the first example are confirmed.
The Hamiltonian Lanczos method and the Arnoldi method perform alike just as the HEKS algorithm and the EKS method. For the functions $\exp$ and  $\cos$  the first two methods converge faster than the latter two. The use of the extended Krylov subspace does not result in faster convergence.
But for the matrix sign function, the two method based on the extended Krylov subspace perform much better.

\section{Concluding Remarks}
The HEKS algorithm for computing a $J$-orthogonal basis of the extended Krylov subspace $\mathcal{K}_{r,s}$ \eqref{eq:extKrylov} has been presented. Unlike the EKSM for generating an orthogonal basis of  $\mathcal{K}_{r,r}$ it allows for short recurrences. The convergence analysis provide in \cite{KniS10} does not apply here as  the field of values of a Hamiltonian matrix does not (strictly) lie in the right half-plane.
Numerical experiments suggest that it may be useful to employ the HEKS algorithm for the approximation of the action of $f(H)$ on a vector $u$ for Hamiltonian matrices $H.$ The performance of the HEKS algorithm is similar to that of EKSM, but HEKS guarantees the structure-preserving projection of the Hamiltonian matrix which may be relevant for some applications.


\appendix \section{Derivation of the HEKS Algorithm}\label{sec:derivation}
This section is devoted to deriving short recurrences for the HEKS algorithm. We will follow the idea sketched in Section \ref{sec:idea}. First $S_1 \in \mathbb{R}^{2n \times 2}$ is constructed such that $S_1^TJ_nS_1 = J_1$ and the columns of $S_1$ span the same subspace as $\mathcal{K}_{2}(H,u_1)$ (that is, $\operatorname{range}\{S_{1}\}=\mathcal{K}_{2}(H,u_1)$). Here $H \in \mathbb{R}^{2n \times 2n}$  is the Hamiltonian matrix under consideration and $u_1\in \mathbb{R}^{2n}$ a given vector with $\|u_1\|_2=1$. Next $S_{2k}\in \mathbb{R}^{2n \times 4k}$ is constructed by extending $S_{2k-1}$ by two columns such that $S_{2k}^TJ_nS_{2k} = J_{2k}$ and  $\operatorname{range}\{S_{2k}\}=\mathcal{K}_{2k}(H,u_1)+ \mathcal{K}_{2k}(H^{-1},H^{-1}u_1).$ Finally, $S_{2k+1}\in \mathbb{R}^{2n \times 4k+2}$ is constructed by extending $S_{2k}$ by two columns such that $S_{2k+1}^TJ_nS_{2k+1} = J_{2k+1}$ and  $\operatorname{range}\{S_{2k+1}\}= \mathcal{K}_{2k+1}(H,u_1)+ \mathcal{K}_{2k}(H^{-1},H^{-1}u_1).$ In doing so, we will provide a proof that the projected matrices
$H_{2k}$ and $H_{2k+1}$ as well as  $J_{2k}^TS_{2k}^TJ_nH^{-1}S_{2k}$ and $J_{2k+1}^TS_{2k+1}^TJ_nH^{-1}S_{2k+1}$ are of the above given forms \eqref{eq_struct1}, \eqref{eq_myHinv} and \eqref{eq_myHinv2}, resp.. In particular, we will prove Theorem \ref{theo1}.
The assumption in Theorem \ref{theo1} that no breakdown occurs in particular implies that in the following all assumptions on nonzero parameters must hold.

\subsection{Step 1: $\operatorname{range}\{S_1\} = \mathcal{K}_{2}(H,u_1)$}
As $u_1$ satisfies $\|u_1\|_2=1,$ there is nothing to do with the first vector in $ \mathcal{K}_{2}(H,u_1)$. The second vector $ Hu_1$ needs to $J$-orthogonalized against $u_1.$ This is achieved by
\begin{equation}\label{eq_v1}
v_1 = Hu_1/u_1^TJ_nHu_1 = Hu_1/\vartheta_1
\end{equation}
assuming that $\vartheta_1 \neq 0.$
 Thus, the matrix $S_1 = [u_1\mid v_1]$ has $J$-orthogonal columns by construction
\begin{align*}
S_1^TJ_nS_1 = \begin{bmatrix}
u_1^TJ_nu_1 & u_1^TJ_nv_1\\ v_1^TJ_nu_1 & v_1^TJ_nv_1
\end{bmatrix}
= \begin{bmatrix}0 & 1\\ -1 & 0\end{bmatrix}
\end{align*}
as any vector is $J$-orthogonal to itself,
$u_1^TJ_nv_1 = u_1^TJ_n Hu_1 / u_1^TJ_nHu_1  = 1$ and $v_1^TJ_nu_1 = (u_1^TJ_n^Tv_1)^T= -(u_1^TJ_nv_1)^T.$
\subsubsection{The projected matrix $H_1 = J_1^TS_1^TJ_nHS_1$}
We will prove that
\begin{align}\label{eq_projS1}
H_1 = J_1^TS_1^TJ_nHS_1 = \begin{bmatrix}
-v_1^TJ_nHu_1 & -v_1^TJ_nHv_1\\
u_1^TJ_nHu_1 & u_1^TJ_nHv_1
\end{bmatrix}
= \begin{bmatrix}
0 & \alpha_1 \\ \vartheta_1& 0
\end{bmatrix}
\end{align}
holds. Due to \eqref{eq_v1} we have $Hu_1 = \vartheta_1 v_1$ and thus
\[
v_1^TJ_nHu_1 = \vartheta_1 \cdot v_1^TJ_nv_1 =0
\]
 as any vector is $J$-orthogonal to itself.
The zero in position (2,2) follows from the zero in position $(1,1)$ as $H$ as well as $H_1$ is Hamiltonian (or by noting that
$ 0 =v_1^TJ_nHu_1 = (v_1^TJ_nH^Tu_1)^T = u_1^T(J_nH)^Tv_1 = u_1^TJ_nHv_1 = 0$).
\subsubsection{The projected matrix $ J_1^TS_1^TJ_nH^{-1}S_1$}
Making use of the fact that $H^TJ_nH^{-1} = -J_n$ as $H$ is Hamiltonian ($(J_nH)^T=-H^TJ_n=J_nH$),
we have due to \eqref{eq_v1}
\begin{align*}
    \vartheta_1 \cdot v_1^TJ_nH^{-1}u_1 = (Hu_1)^TJ_nH^{-1}u_1
    = u_1^T H^TJ_nH^{-1}u_1
    = -u_1^TJ_nu_1 = 0.
\end{align*}
This implies $u_1^TJ_nH^{-1}v_1 =0.$ Moreover, using  \eqref{eq_v1} again
\begin{align*}
v_1^TJ_nH^{-1}v_1 &= (Hu_1)^TJ_nH^{-1} (Hu_1)/\vartheta_1^2 = u_1^TH^TJ_nu_1/\vartheta_1^2 = -u_1^TJ_nHu_1/\vartheta_1^2 =-1/\vartheta_1.
\end{align*}
Thus
\begin{equation}\label{eq_hinv1}
J_1^TS_1^TJ_nH^{-1}S_1 = \begin{bmatrix}
-v_1^TJ_nH^{-1}u_1 & -v_1^TJ_nH^{-1}v_1\\
u_1^TJ_nH^{-1}u_1 & u_1^TJ_nH^{-1}v_1
\end{bmatrix}
= \begin{bmatrix}
0 & 1/\vartheta_1 \\ f_{11}& 0
\end{bmatrix}.
\end{equation}
\subsection{Step 2:  $\operatorname{range}\{S_2\} = \mathcal{K}_{2}(H,u_1)+ \mathcal{K}_{2}(H^{-1},H^{-1}u_1)$}
Now the first vector from the Krylov subspace $\mathcal{K}_{r}(H^{-1},H^{-1}u_1)$ is added to the symplectic basis by $J$-orthogonalization $H^{-1}u_1$ against $u_1$ and $v_1.$ This is achieved by computing
\begin{align*}
    w_x = (I-S_1J_1^TS_1^TJ_n)H^{-1}u_1,
\end{align*}
and normalizing $w_x$ to length 1,  $x_1 = w_x/ \|w_x\|_2$.
Next the second vector from $\mathcal{K}_{r}(H^{-1},H^{-1}u_1)$ needs to be added to the symplectic basis.
This can be accomplished by $J$-orthogonalizing $H^{-1}x_1$  against $u_1$ and $v_1$
\begin{align*}
    w_y &=(I-S_1J_1^TS_1^TJ_n)H^{-1}x_1
\end{align*}
and making sure that $w_y$ is $J$-orthogonal against $x_1$ as well, $y_1 = w_y/w_y^TJx_1.$ Here we assume that $\|w_x\|_2 \neq 0$ as well as $w_y^TJx_1 \neq 0.$

Collect the vectors into a matrix $S_2 =[y_1~u_1\mid x_1~v_1] \in \mathbb{R}^{2n \times 4}.$ By construction the columns of $S_2$ are  $J$-orthogonal, that is
\begin{align}\label{eq_app1}
S_2^TJ_nS_2=
J_2
\end{align}
and
\[
\operatorname{range}\{S_2\} = \mathcal{K}_{2}(H,u_1)+ \mathcal{K}_{2}(H^{-1},H^{-1}u_1).
\]

Let us take a closer look at $w_x$ and $w_y.$ Making use of \eqref{eq_hinv1} we have
\begin{align*}
    w_x &= H^{-1}u_1 -[v_1~-u_1]\begin{bmatrix}
    u_1^TJ_nH^{-1}u_1\\v_1^TJ_nH^{-1}u_1
    \end{bmatrix}
= H^{-1}u_1 -[v_1~-u_1]\begin{bmatrix}
    f_{11}\\0
    \end{bmatrix}
=  H^{-1}u_1 -f_{11} v_1. 
\end{align*}
Hence, with $\psi_1 = \|w_x\|_2$
we have
\begin{equation}\label{eq_x1}
 x_1 = \left(H^{-1}u_1 -f_{11} v_1\right )/\psi_1,
\end{equation}
where,  as already stated above, $\psi_1 \neq 0$ is assumed.

Next we turn our attention to $w_y.$ We will make use of the fact that $H^{-1}$ is Hamiltonian  ($J_nH^{-1}=-H^{-T}J_n$) and $S_2^TJ_nS_2=J_2.$
With \eqref{eq_x1} we see
\begin{align}
 u_1^TJ_nH^{-1}x_1 &= -(H^{-1}u_1)^TJ_n x_1
= -(\psi_1 x_1+f_{11}v_1)^TJ_nx_1=0.
\label{eq_1}
\end{align}
Similarly, it follows  with \eqref{eq_v1} that
\begin{align}\label{eq_2}
v_1^TJ_nH^{-1}x_1 &= -(H^{-1}v_1)^TJ_nx_1 = u_1^TJ_nx_1/\vartheta_1 = 0.
\end{align}
Hence,
\begin{align*}
    w_y &= H^{-1}x_1 - [v_1~-u_1]\begin{bmatrix}
    u_1^TJ_nH^{-1}x_1\\v_1^TJ_nH^{-1}x_1
    \end{bmatrix}
= H^{-1}x_1 - [v_1~-u_1]\begin{bmatrix}
    0\\0
    \end{bmatrix}
 = H^{-1}x_1,
\end{align*}
and
\[
y_1 = H^{-1}x_1/\xi_1,
\]
where we assume that
\begin{equation}\label{eq_xi1}
\xi_1 = (H^{-1}x_1)^TJ_nx_1 = x_1^TH^{-T}J_nx_1\neq 0.
\end{equation}
Observe that
\begin{align}\label{eq_xi}
\begin{split}
    \delta _1 = y_1^TJ_nHy_1 &
    =\frac{1}{\xi_1^2}x_1^TH^{-T}J_nHH^{-1}x_1
    =\frac{1}{\xi_1^2}x_1^TH^{-T}J_nx_1= \frac{1}{\xi_1}.
\end{split}
\end{align}
Thus
\begin{equation}\label{eq_y1}
y_1 = H^{-1}x_1/\xi_1 = \delta_1 H^{-1}x_1.
\end{equation}
\subsubsection{The projected matrix $H_2 = J_2^TS_2^TJ_nHS_2$}
We will see that the zero structure of  $H_2 = J_2^TS_2^TJ_nHS_2$ is given as follows:
\begin{align}\label{eq_H2}
H_2 &=
\left[\begin{array}{cc||cc}
-x_1^TJ_nHy_1 & -x_1^TJ_nHu_1 & -x_1^TJ_nHx_1 & -x_1^TJ_nHv_1\\
-v_1^TJ_nHy_1 & -{\color{blue}v_1^TJ_nHu_1} & -v_1^TJ_nHx_1 & -{\color{blue}v_1^TJ_nHv_1}\\\hline\hline
y_1^TJ_nHy_1 & y_1^TJ_nHu_1 & y_1^TJ_nHx_1 & y_1^TJ_nHv_1\\
u_1^TJ_nHy_1 & {\color{blue}u_1^TJ_nHu_1} & u_1^TJ_nHx_1 & {\color{blue}u_1^TJ_nHv_1}
\end{array}\right]
=
\left[\begin{array}{cc||cc}
0 & 0 & \lambda_1 & \gamma_1\\
0 &{\color{blue} 0 }&\gamma_1 & {\color{blue}\alpha_1}\\ \hline\hline
\delta_1 & 0 &0 & 0\\
0 & {\color{blue}\vartheta_1 }& 0 &{\color{blue}0}
\end{array}\right].
\end{align}
The entries at the positions $(2,2), (4,2), (2,4)$ and $(4,4)$ (denoted in blue in \eqref{eq_H2}) are the same as in \eqref{eq_projS1}. Due to $H$ and thus $H_2$ being Hamiltonian, we only need to prove the zero entries at the positions $(1,1), (1,2), (2,1)$ and $(3,2)$,
the other zeros in \eqref{eq_H2} follow immediately.

 Due to \eqref{eq_v1}  we have $Hu_1 = \vartheta_1 v_1$. Thus,
$x_1^TJ_nHu_1 = \vartheta_1 \cdot x_1^TJ_nv_1 = 0$ and $y_1^TJ_nHu_1 = \vartheta_1 \cdot y_1^TJ_nv_1 = 0$ due to \eqref{eq_app1}. This gives the zero entries in the positions $(1,2)$ and $(3,2)$.

Due to \eqref{eq_y1} it follows with \eqref{eq_app1}  for the entry $(2,1)$
\begin{align*}
   v_1^T J_n H y_1 &=\delta_1  v_1^T J_n H H^{-1}x_1
    = \delta_1 v_1^T J_n x_1 =0.
\end{align*}
Moreover,  in a similar way for the entry $(1,1)$ we have
\begin{align*}
   x_1^T J_n H y_1 &=\delta_1  x_1^T J_n H H^{-1}x_1
    = \delta_1 x_1^T J_n x_1 =0.
\end{align*}
Hence, \eqref{eq_H2} holds.
\subsubsection{The projected matrix $ J_2^TS_2^TJ_nH^{-1}S_2$}
Some of the entries in $\tilde H_2 = J_2^TS_2^TJ_nH^{-1}S_2$ (denoted in blue) are already known from \eqref{eq_hinv1},
\begin{align}
\tilde H_2 &=
\left[\begin{array}{cc|cc}
-x_1^TJ_nH^{-1}y_1 & -x_1^TJ_nH^{-1}u_1 & -x_1^TJ_nH^{-1}x_1 & -x_1^TJ_nH^{-1}v_1\\
-v_1^TJ_nH^{-1}y_1 & -{\color{blue}v_1^TJ_nH^{-1}u_1} & -v_1^TJ_nH^{-1}x_1 & -{\color{blue}v_1^TJ_nH^{-1}v_1}\\\hline
y_1^TJ_nH^{-1}y_1 & y_1^TJ_nH^{-1}u_1 & y_1^TJ_nH^{-1}x_1 & y_1^TJ_nH^{-1}v_1\\
u_1^TJ_nH^{-1}y_1 & {\color{blue}u_1^TJ_nH^{-1}u_1} & u_1^TJ_nH^{-1}x_1 & {\color{blue}u_1^TJ_nH^{-1}v_1}
\end{array}\right] \nonumber\\ &
=
\left[\begin{array}{cc|cc}
0 &0 & 1/\delta_1 & 0\\
 0&  {\color{blue}0} & 0&  {\color{blue}1/\vartheta_1}\\ \hline
e_{11} & g_{11} & 0& 0\\
g_{11} &  {\color{blue} f_{11}}& 0 & {\color{blue}0}
\end{array}\right].\label{eq_hinv2}
\end{align}
The entry in position $(1,3)$ follows from \eqref{eq_xi1} and \eqref{eq_xi}, while the zero entries in the positions $(1,2)$, $(1,4)$, $(2,3)$ and $(4,3)$ have already been proven in \eqref{eq_1} and \eqref{eq_2}.

It remains to consider the entries at the positions $(3,3)$ and $(3,4)$. Using \eqref{eq_y1} and \eqref{eq_v1} leads to
\begin{align*}
y_1^TJ_nH^{-1}x_1 &= \delta_1 (H^{-1}x_1)^TJ_n(H^{-1}x_1) = 0,\\
 y_1^TJ_nH^{-1}v_1 &= y_1^TJ_nu_1/\vartheta_1 = 0.
\end{align*}
Hence, \eqref{eq_hinv2} holds.
\subsection{Step 3:  $\operatorname{range}\{S_3\} = \mathcal{K}_{4}(H,u_1)+ \mathcal{K}_{2}(H^{-1},H^{-1}u_1)$}

In this step the next two vectors $H^2u_1$ and $H^3u_1$ from $\mathcal{K}_4(H,u_1)$ are added to the symplectic basis. We start by $J$-orthogonalizing $ Hv_1$ against the columns of $S_2$
\begin{align}
    w_u &=(I-S_2J_2^TS_2^TJ_n)Hv_1
    = Hv_1 -[x_1~~v_1~~-y_1~~-u_1]\begin{bmatrix}
    y_1^TJ_nHv_1\\
    u_1^TJ_nHv_1\\
    x_1^TJ_nHv_1\\
    v_1^TJ_nHv_1
    \end{bmatrix} \nonumber \\&
= Hv_1 -[x_1~~v_1~~-y_1~~-u_1]\begin{bmatrix}
   0\\
   0\\
   -\gamma_1\\
    -\alpha_1
    \end{bmatrix}
= Hv_1-\gamma_1y_1 - \alpha_1u_1\label{eq_wu}
\end{align}
where \eqref{eq_H2} gives that the first two entries of the last vector are zero.
Normalizing $w_u$ to length $1$  gives $u_2$
\begin{equation}\label{eq_u2}
u_2 = w_u/\chi_2,
\end{equation}
where it is assumed that $\chi_2 = \|w_u\|_2\neq 0.$

This step is finalized by $J$-orthogonalizing $Hu_2$ against the columns of $S_2$:
\begin{align*}
    w_v &=(I-S_2J_2^TS_2^TJ_n)Hu_2
    = Hu_2 -[x_1~~v_1~~-y_1~~-u_1]\begin{bmatrix}
    y_1^TJ_nHu_2\\
    u_1^TJ_nHu_2\\
    x_1^TJ_nHu_2\\
    v_1^TJ_nHu_2
    \end{bmatrix}.
\end{align*}
All entries of the last vector are zero. The first two zeros follow as $H^{-T}J_nH=-J_n$  with \eqref{eq_y1} and \eqref{eq_v1}:
\begin{align*}
   y_1^TJ_nHu_2 /\delta_1 &=  (H^{-1}x_1)^TJ_nHu_2 =  -x_1^TJ_nu_2 = 0,\\
  u_1^TJ_nHu_2/\vartheta_1  &= (H^{-1}v_1)^TJ_nHu_2 =-v_1^TJ_nu_2 = 0
\end{align*}
by construction of $u_2.$ The last zero follows as $H$ is Hamiltonian with \eqref{eq_wu},
\begin{align*}
  v_1^TJ_nHu_2 &=v_1^T(J_nH)^Tu_2 = -(Hv_1)^TJ_nu_2
 =-(\chi_2u_2+\gamma_1y_1+\alpha_1u_1)^TJ_nu_2 = 0
 \end{align*}
again due to the construction of $u_2.$  With this and \eqref{eq_x1} we have for the next to  last entry
\begin{align*}
  \psi_1 \cdot x_1^TJ_nHu_2 &=(H^{-1}u_1+f_{11} v_1)^TJ_nHu_2
= -u_1^TJ_nu_2+f_{11} v_1^TJ_nHu_2 = 0.
\end{align*}
Thus the expression for $w_v$ simplifies to
\[
w_v = Hu_2.
\]
Normalizing $w_v$ by $\vartheta_2 = u_2^TJ_nHu_2$ to make sure it is $J$-orthogonal to $u_2$ as well yields
\[v_2 = Hu_2/\vartheta_2.\]
Let
\[
S_3 = [y_1~u_1~u_2\mid x_1~v_1~v_2] \in \mathbb{R}^{2n \times 6}.
\]
Then by construction
\begin{equation}\label{eq_S3symp}
S_3^TJ_nS_3 = J_3
\end{equation}
and
\[
\operatorname{range}\{S_3\} = \mathcal{K}_{4}(H,u_1)+ \mathcal{K}_{2}(H^{-1},H^{-1}u_1).
\]
\subsubsection{The projected matrix $H_3 = J_3^TS_3^TJ_nHS_3$}
Some of the entries (denoted in blue) in $H_3 = J_3^TS_3^TJ_nHS_3$ are already known from \eqref{eq_H2}
\begin{align}
H_3 &=
{\small \left[\begin{array}{c|cc||c|cc}
{\color{blue}0} & {\color{blue}0} & -x_1^TJ_nHu_2 & {\color{blue}\lambda_1} & {\color{blue}\gamma_1} & -x_1^TJ_nHv_2\\ \hline
{\color{blue}0} & {\color{blue}0} &-v_1^TJ_nHu_2 & {\color{blue}\gamma_1} & {\color{blue}\alpha_1} & -v_1^TJ_nHv_2\\
-v_2^TJ_nHy_1 & -v_2^TJ_nHu_1 & -v_2^TJ_nHu_2 & -v_2^TJ_nHx_1 & -v_2^TJ_nHv_1 & -v_2^TJ_nHv_2\\
\hline\hline
{\color{blue}\delta_1} &{\color{blue} 0} &y_1^TJ_nHu_2 & {\color{blue}0} & {\color{blue}0} & y_1^TJ_nHv_2\\\hline
{\color{blue}0} & {\color{blue}\vartheta_1} & u_1^TJ_nHu_2 &{\color{blue}0} &{\color{blue}0} & u_1^TJ_nHv_2\\
u_2^TJ_nHy_1 & u_2^TJ_nHu_1 &u_2^TJ_nHu_2 &u_2^TJ_nHx_1 &u_2^TJ_nHv_1 &u_2^TJ_nHv_2
\end{array}\right] }\nonumber\\
&=
\left[\begin{array}{c|cc||c|cc}
{\color{blue}0} & {\color{blue}0} & 0 &{\color{blue}\lambda_1} & {\color{blue}\gamma_1} & \mu_2\\\hline
{\color{blue}0} & {\color{blue}0} &0 & {\color{blue}\gamma_1} & {\color{blue}\alpha_1} & \beta_2\\
0 & 0 & 0 & \mu_2 & \beta_2 & \alpha_2\\
\hline\hline
{\color{blue}\delta_1} & {\color{blue}0} &0 & {\color{blue}0} & {\color{blue}0} & 0\\\hline
{\color{blue}0} & {\color{blue}\vartheta_1} & 0 &{\color{blue}0} &{\color{blue}0} & 0\\
0 & 0 &\vartheta_2 &0 &0 &0
\end{array}\right].\label{eq_H3}
\end{align}
The zeros in the third column (and hence the zeros in the last row) follow with $Hu_2 = \vartheta_2 v_2$ due to \eqref{eq_S3symp}.
Moreover, we have with \eqref{eq_y1} and \eqref{eq_v1}
\begin{align*}
v_2^TJ_nHy_1 &= \delta_1 v_2^TJ_n x_1 = 0,\\
v_2^TJ_nHu_1 &= \vartheta_1 v_2^TJ_nv_1 =0
\end{align*}
making again use of \eqref{eq_S3symp}. Hence, \eqref{eq_H3} holds.
\subsubsection{The projected matrix $ J_3^TS_3^TJ_nH^{-1}S_3$}
Some of the entries in $\tilde H_3 = J_3^TS_3^TJ_nH^{-1}S_3$ (denoted in blue)  are already known from \eqref{eq_hinv2}
\begin{align}
\tilde H_3 &=
{\footnotesize \left[\begin{array}{ccc|ccc}
{\color{blue}0} & {\color{blue}0} & -x_1^TJ_nH^{-1}u_2 & {\color{blue}1/\delta_1} &  {\color{blue}0}&-x_1^TJ_nH^{-1}v_2\\
{\color{blue}0} &  {\color{blue}0} &-v_1^TJ_nH^{-1}u_2 &  {\color{blue}0} &  {\color{blue}1/\vartheta_1} & -v_1^TJ_nH^{-1}v_2\\
-v_2^TJ_nH^{-1}y_1 & -v_2^TJ_nH^{-1}u_1 & -v_2^TJ_nH^{-1}u_2 & -v_2^TJ_nH^{-1}x_1 & -v_2^TJ_nH^{-1}v_1 & -v_2^TJ_nH^{-1}v_2\\
\hline
 {\color{blue}e_{11}} &  {\color{blue}g_{11}} &y_1^TJ_nH^{-1}u_2 &  {\color{blue}0} &  {\color{blue}0} & y_1^TJ_nH^{-1}v_2\\
 {\color{blue}g_{11}} &  {\color{blue}f_{11}} & u_1^TJ_nH^{-1}u_2 &0 &0 & u_1^TJ_nH^{-1}v_2\\
u_2^TJ_nH^{-1}y_1 & u_2^TJ_nH^{-1}u_1 &u_2^TJ_nH^{-1}u_2 &u_2^TJ_nH^{-1}x_1 &u_2^TJ_nH^{-1}v_1 &u_2^TJ_nH^{-1}v_2
\end{array}\right] }\nonumber\\
&=
\left[\begin{array}{ccc|ccc}
 {\color{blue}0 }&  {\color{blue}0} & 0 &  {\color{blue}1/\delta_1} &  {\color{blue}0} & 0\\
 {\color{blue}0} &  {\color{blue}0} &0 &  {\color{blue}0} &  {\color{blue}1/\vartheta_1} & 0\\
0 & 0 & 0 & 0 & 0 & 1/\vartheta_2\\
\hline
 {\color{blue}e_{11}} &  {\color{blue}g_{11}} &g_{12} &  {\color{blue}0} &  {\color{blue}0} & 0\\
 {\color{blue}g_{11}} &  {\color{blue}f_{11}} & 0& {\color{blue}0} & {\color{blue}0} & 0\\
g_{12} & 0 &f_{22} &0 &0 &0
\end{array}\right].\label{eq_hinv3}
\end{align}
It remains to show that the five entries $v_2^TJ_nH^{-1}z$ for $z =x_1, v_1,  y_1, u_1, u_2$ as well as the three entries $z^TJ_nH^{-1}u_2$ for $z = v_1,x_1,u_1$ are zero. Moreover, we need to show that $-v_2^TJ_nH^{-1}v_2=1/\vartheta_2.$

Most of these relations follow from $H^TJ_nH^{-T}=-J_n$ and due to $Hu_j =\vartheta_j v_j, j = 1,2.$
Making use of \eqref{eq_S3symp} in the last equality of each equation
 we have
\begin{align*}
\vartheta_1 \cdot v_1^TJ_nH^{-1}u_2 &= u_1^TH^TJ_nH^{-1}u_2 = -u_1^TJ_nu_2 = 0,\\
\vartheta_2 \cdot v_2^TJ_nH^{-1}v_2 &= u_2^TH^TJ_nH^{-1}v_2 = -u_2^TJ_nv_2 = -1,
\end{align*}
and for $z =x_1, v_1,  y_1, u_1, u_2$
\begin{align*}
\vartheta_2 \cdot v_2^TJ_nH^{-1}z &= u_2^TH^TJ_nH^{-1}z = -u_2^TJ_nz = 0.
\end{align*}
Thus, $- v_2^TJ_nH^{-1}v_2 = 1/\vartheta_2$ and the five entries in the $(1,1)$ (and the $(2,2)$) block of $\tilde H_3$ are zero.

The derivation of the final two zero entries needs a slightly more involved derivation.
Due to \eqref{eq_u2}, \eqref{eq_S3symp} and \eqref{eq_hinv2} we have
\begin{align*}
\chi_2 \cdot x_1^TJ_nH^{-1}u_2 &= x_1^TJ_nH^{-1} (Hv_1-\gamma_1y_1-\alpha_1u_1)\\
&= x_1^TJ_nv_1-\gamma_1x_1^TJ_nH^{-1}y_1-\alpha_1x_1^TJ_nH^{-1}u_1=0,
\end{align*}
while due to $H^{-1}$ being Hamiltonian, \eqref{eq_x1} and \eqref{eq_S3symp} we get
\[
 u_1^TJ_nH^{-1}u_2 = (H^{-1}u_1)^TJ_nu_2
 = (\psi_1x_1+f_{11}v_1)^TJ_nu_2
=0.
\]
Hence, \eqref{eq_hinv3} holds.
\subsection{Step 4:  $\operatorname{range}\{S_4\} = \mathcal{K}_{4}(H,u_1)+ \mathcal{K}_{4}(H^{-1},H^{-1}u_1)$}

In this step the next two vectors $H^{-3}u_1$ and $H^{-4}u_1$ from $\mathcal{K}_4(H^{-1},H^{-1}u_1)$ are added to the symplectic basis. We start by $J$-orthogonalizing $H^{-1}y_1$ against the columns of $S_3$
\begin{align*}
    w_x &=(I-S_3J_3^TS_3^TJ_n)H^{-1}y_1
    = H^{-1}y_1 -[x_1~V_2\mid -y_1~-U_2]\begin{bmatrix}
    y_1^TJ_nH^{-1}y_1\\
    U_2^TJ_nH^{-1}y_1\\
    x_1^TJ_nH^{-1}y_1\\
    V_2^TJ_nH^{-1}y_1
    \end{bmatrix} \nonumber\\
    &= H^{-1}y_1 -[x_1~V_2\mid -y_1~-U_2]\begin{bmatrix}
    e_{11}\\g_{11}\\g_{12}\\0\\0\\0
    \end{bmatrix} 
= H^{-1}y_1 -e_{11}x_1-g_{11}v_1-g_{12}v_2 
\end{align*}
due to \eqref{eq_hinv3}.

Normalizing $w_x$ to length $1$  gives
\begin{equation}\label{eq_x2}
x_2 = w_x/\psi_2,
\end{equation}
where we assume that $\psi_2 =\|w_x\|_2\neq0.$

This step is finalized by $J$-orthogonalizing $H^{-1}x_2$ against the columns of $S_3$
\begin{align}\label{eq_xyz}
    w_y&=(I-S_3J_3^TS_3^TJ_n)H^{-1}x_2
    = H^{-1}x_2 -[x_1~V_2\mid -y_1~-U_2]\begin{bmatrix}
    y_1^TJ_nH^{-1}x_2\\
    U_2^TJ_nH^{-1}x_2\\
    x_1^TJ_nH^{-1}x_2\\
    V_2^TJ_nH^{-1}x_2
    \end{bmatrix}.
\end{align}
All entries $z^TJ_nH^{-1}x_2 =-(H^{-1}z)^TJ_nx_2$ in the last vector are zero.
As $Hu_i = \vartheta_i v_i,$ $ i = 1,2,$ for the last two entries we have
\begin{align*}
\vartheta_1 \cdot (H^{-1}v_1)^TJ_nx_2 &= u_1^TJ_nx_2 =0,\\
\vartheta_2 \cdot (H^{-1}v_2)^TJ_nx_2 &= u_2^TJ_nx_2 =0
\end{align*}
by construction of $x_2$  as a vector $J$-orthogonal to all columns of $S_3.$
Next, we use \eqref{eq_x2}, \eqref{eq_y1} and \eqref{eq_x1} to see
\begin{align*}
(H^{-1}y_1)^TJ_nx_2 &= (\psi_2x_2-e_{11}x_1-g_{11}v_1-g_{12}v_2)^TJ_nx_2 = 0,\\
(H^{-1}x_1)^TJ_nx_2 &= \xi_1y_1^TJ_nx_2 = 0,\\
(H^{-1}u_1)^TJ_nx_2 &= (\psi_1x_1-f_{11}v_1)^TJ_nx_2 =0
\end{align*}
again by construction of $x_2$  as a vector $J$-orthogonal to all columns of $S_3.$ With this and \eqref{eq_u2}, it follows that
\begin{align*}
\chi_2 \cdot (H^{-1}u_2)^TJ_nx_2 &= (v_1-\gamma_1H^{-1}y_1-\alpha_1H^{-1}u_1)^TJ_nx_2 =0.
\end{align*}
Thus,
\[
y_2 = H^{-1}x_2/\xi_2,
\]
where we assume that
\begin{equation}\label{eq_xi2}
\xi_2 = (H^{-1}x_2)^TJ_nx_2 = x_2^TH^{-T}J_nx_2\neq 0.
\end{equation}
With the same argument as in \eqref{eq_xi} we see that
\[
    \delta _2 = y_2^TJ_nHy_2 = \frac{1}{\xi_2^2}(H^{-1}x_2)^TJ_nHH^{-1}x_2 = \frac{1}{\xi_2}.
\]
Thus
\begin{equation}\label{eq_y2}
y_2 = H^{-1}x_2/\xi_2 =\delta_2 H^{-1}x_2.
\end{equation}
Let
\[
S_4 = [y_2~y_1~u_1~u_2\mid x_2~x_1~v_1~v_2] \in \mathbb{R}^{2n \times 8}.
\]
Then by construction
\begin{equation}\label{eq_S4symp}
S_4^TJ_nS_4 = J_4
\end{equation}
and
\[
\operatorname{range}\{S_4\} = \mathcal{K}_{4}(H,u_1)+ \mathcal{K}_{4}(H^{-1},H^{-1}u_1).
\]
\subsubsection{The projected matrix $H_4 = J_4^TS_4^TJ_nHS_4$}
Some of the entries in $H_4 = J_4^TS_4^TJ_nHS_4$ (denoted in blue) are already known from \eqref{eq_H3}
\[
H_4 =
{\footnotesize \left[\begin{array}{cc|cc||cc|cc}
-x_2^TJ_nHy_2 & -x_2^TJ_nHy_1 & -x_2^TJ_nHu_1  & -x_2^TJ_nHu_2 & -x_2^TJ_nHx_2 & -x_2^TJ_nHx_1 & -x_2^TJ_nHv_1 & -x_2^TJ_nHv_2 \\
-x_1^TJ_nHy_2 & {\color{blue}0} & {\color{blue}0} &{\color{blue}0} &-x_1^TJ_nHx_2  & {\color{blue} \lambda_1} & {\color{blue}\gamma_1} & {\color{blue}\mu_2}\\ \hline
-v_1^TJ_nHy_2 & {\color{blue}0} & {\color{blue}0} &{\color{blue}0} & -v_1^TJ_nHx_2 &{\color{blue}\gamma_1} & {\color{blue}\alpha_1} &{\color{blue}\beta_2}\\
-v_2^TJ_nHy_2 &{\color{blue}0} & {\color{blue}0} &{\color{blue}0} & -v_2^TJ_nHx_2 & {\color{blue}\mu_2} & {\color{blue}\beta_2} & {\color{blue}\alpha_2}\\
\hline\hline
y_2^TJ_nHy_2 & y_2^TJ_nHy_1 & y_2^TJ_nHu_1  & y_2^TJ_nHu_2 & y_2^TJ_nHx_2 & y_2^TJ_nHx_1 & y_2^TJ_nHv_1 & y_2^TJ_nHv_2 \\
y_1^TJ_nHy_2&{\color{blue}\delta_1} & {\color{blue}0} &{\color{blue}0} &y_1^TJ_nHx_2 &{\color{blue}0} & {\color{blue}0} & {\color{blue}0}\\\hline
u_1^TJ_nHy_2&{\color{blue}0} & {\color{blue}\vartheta_1} & {\color{blue}0} &u_1^TJ_nHx_2&{\color{blue}0} &{\color{blue}0} & {\color{blue}0}\\
u_2^TJ_nHy_2&{\color{blue}0} & {\color{blue}0}& {\color{blue}\vartheta_2} &u_2^TJ_nHx_2 &{\color{blue}0} &{\color{blue}0} &{\color{blue}0}
\end{array}\right] }
\]
We will show that
\begin{align}
H_4=
\left[\begin{array}{cc|cc||cc|cc}
0 & 0 & 0 & 0 & \lambda_2 & 0 & 0 & \gamma_2\\
0& {\color{blue}0} & {\color{blue}0} & {\color{blue}0} & 0& {\color{blue}\lambda_1} & {\color{blue}\gamma_1} & {\color{blue}\mu_2}\\\hline
0&  {\color{blue}0} & {\color{blue}0} &{\color{blue}0} & 0& {\color{blue}\gamma_1} & {\color{blue}\alpha_1} & {\color{blue}\beta_2}\\
0& {\color{blue}0} & {\color{blue}0} & {\color{blue}0} & \gamma_2& {\color{blue}\mu_2} & {\color{blue}\beta_2} & {\color{blue}\alpha_2}\\
\hline\hline
\delta_2 &0 & 0 & 0&0&0&0&0\\
0 &{\color{blue}\delta_1} & {\color{blue}0} &{\color{blue}0} & 0 & {\color{blue}0} & {\color{blue}0}&{\color{blue}0}\\\hline
0 & {\color{blue}0} & {\color{blue}\vartheta_1} & {\color{blue}0} &0 &{\color{blue}0} & {\color{blue}0}&{\color{blue}0}\\
0 &{\color{blue} 0} & {\color{blue}0} &{\color{blue}\vartheta_2} &0 &{\color{blue}0} &{\color{blue}0}&{\color{blue}0}
\end{array}\right].\label{eq_H4}
\end{align}
Let us consider the entries in the first column of \eqref{eq_H4}. We make use of \eqref{eq_y2} and obtain
\[
z^TJ_nHy_2 = \delta_2 z^TJ_nx_2.
\]
For $z =x_1, x_2, v_1, v_2, y_1,u_1, u_2$ we have $z^TJ_nx_2=0$ due to  \eqref{eq_S4symp}, while,
$y_2^TJ_nx_2 =1.$ Thus $y_2^TJ_nHy_2 = \delta_2.$ Moreover, the other 7 entries in the first column are zero. This implies that the other  7 entries in the fifth row are zero as well.

For the entries $x_2^TJ_nHz$ for $z = y_1, u_1, u_2,$ in the first row we note that
\begin{align*}
\psi_2 \cdot x_2^TJ_nHz &= (H^{-1}y_1- e_{11}x_1-g_{11}v_1-g_{12}v_2)^TJ_nHz
= - y_1^TJ_nz-( e_{11}x_1+g_{11}v_1+g_{12}v_2)^TJ_nHz =0
\end{align*}
due to \eqref{eq_x2}, \eqref{eq_S4symp} and \eqref{eq_H3}. These zeros imply zeros in the positions $(6,5), (7,5)$ and $(8,5).$

It remains to consider the two entries at the positions $(1,6)$ and $(1,7)$. With \eqref{eq_x1} it follows that
\begin{align*}
\psi_1 \cdot x_2^TJ_nHx_1 &= x_2^TJ_nH(H^{-1}u_1-f_{11}v_1)
= x_2^TJ_nu_1 - f_{11}x_2^TJ_nHv_1
= -f_{11}x_2^TJ_nHv_1
\end{align*}
due to \eqref{eq_S4symp}. Thus,
the entry at position $(1,6)$ is zero if and only if the entry at position $(1,7)$ is zero. For the entry at position $(1,7)$ we have with \eqref{eq_u2} and \eqref{eq_wu}
\begin{align*}
x_2^TJ_nHv_1 &= x_2^TJ_n(\chi_2u_2+\gamma_1y_1\alpha_1u_1) = 0
\end{align*}
due to \eqref{eq_S4symp}. Hence, \eqref{eq_H4} holds.
\subsubsection{The projected matrix $J_4^TS_4^TJ_nH^{-1}S_4$}
Some of the entries in $\tilde H_4 = J_4^TS_4^TJ_nH^{-1}S_4$ (denoted in blue) are already known from \eqref{eq_hinv3},
\begin{align}
&
{\tiny \left[\begin{array}{cccc|cccc}
-x_2^TJ_nH^{-1}y_2 & {\color{red}0} & {\color{red}0}  &{\color{red}0} & -x_2^TJ_nH^{-1}x_2 & {\color{red}0} & {\color{red}0} & {\color{red}0} \\
-x_1^TJ_nH^{-1}y_2 & {\color{blue}0} & {\color{blue}0} &{\color{blue}0} &{\color{red}0}  &  {\color{blue}1/\delta_1} &
{\color{blue}0} & {\color{blue}0}\\
-v_1^TJ_nH^{-1}y_2 & {\color{blue}0} & {\color{blue}0} &{\color{blue}0} & {\color{red}0} &{\color{blue}0} & {\color{blue}1/\vartheta_1} &{\color{blue}0}\\
-v_2^TJ_nH^{-1}y_2 &{\color{blue}0} & {\color{blue}0} &{\color{blue}0} & {\color{red}0} & {\color{blue}0} & {\color{blue}0} & {\color{blue}1/\vartheta_2}\\
\hline
y_2^TJ_nH^{-1}y_2 & y_2^TJ_nH^{-1}y_1 & y_2^TJ_nH^{-1}u_1  & y_2^TJ_nH^{-1}u_2 & y_2^TJ_nH^{-1}x_2 & y_2^TJ_nH^{-1}x_1 & y_2^TJ_nH^{-1}v_1 & y_2^TJ_nH^{-1}v_2 \\
y_1^TJ_nH^{-1}y_2&{\color{blue}e_{11}} & {\color{blue}g_{11}} &{\color{blue}g_{12}} &{\color{red}0} &{\color{blue}0} & {\color{blue}0} & {\color{blue}0}\\
u_1^TJ_nH^{-1}y_2&{\color{blue}g_{11}} & {\color{blue}f_{11}} & {\color{blue}0} &{\color{red}0}&{\color{blue}0} &{\color{blue}0} & {\color{blue}0}\\
u_2^TJ_nH^{-1}y_2&{\color{blue}g_{12}} & {\color{blue}0}& {\color{blue}f_{22}} &{\color{red}0} &{\color{blue}0} &{\color{blue}0} &{\color{blue}0}
\end{array}\right] }\nonumber\\
&=
\left[\begin{array}{cc|cc||cc|cc}
0 & {\color{red}0}& {\color{red}0} & {\color{red}0} & 1/\delta_2 & {\color{red}0} & {\color{red}0} & {\color{red}0}\\
0&{\color{blue} 0} & {\color{blue}0} & {\color{blue}0} & {\color{red}0}&  {\color{blue}1/\delta_1} & {\color{blue}0}&{\color{blue}0}\\\hline
0&  {\color{blue}0} &{\color{blue} 0} &{\color{blue}0} & {\color{red}0}& {\color{blue}0} & {\color{blue}1/\vartheta_1} & {\color{blue}0}\\
0& {\color{blue}0} & {\color{blue}0} & {\color{blue}0} & {\color{red}0}&{\color{blue}0} & {\color{blue}0} & {\color{blue}1/\vartheta_2}\\
\hline\hline
e_{22} &e_{12} & g_{12} & g_{22}&0&0&0&0\\
e_{12} &{\color{blue}e_{11}} & {\color{blue}g_{11}} &{\color{blue}g_{12}} & {\color{red}0} & {\color{blue}0} & {\color{blue}0}&{\color{blue}0}\\\hline
g_{12} & {\color{blue}g_{11}}& {\color{blue}f_{11}} & {\color{blue}0} &{\color{red}0} &{\color{blue}0} & {\color{blue}0}&{\color{blue}0}\\
g_{22} & {\color{blue}g_{12}} & {\color{blue}0} &{\color{blue}f_{22}} &{\color{red}0} &{\color{blue}0} &{\color{blue}0}&{\color{blue}0}
\end{array}\right].\label{eq_hinv4}
\end{align}
In addition, most of the ones in the 5th column (denoted in red) and hence in the first row are known from the derivations concerning \eqref{eq_xyz}.

Let us consider the remaining entries in the first row. As $H^{-1}$ is Hamiltonian and due to \eqref{eq_xi2}, we have
\begin{align*}
x_2^TJ_nH^{-1}x_2 &= x_2^T(J_nH^{-1})^Tx_2 = -(H^{-1}x_2)^TJ_nx_2 = -\xi_2 = -1/\delta_2.
\end{align*}
Due to \eqref{eq_y2} and \eqref{eq_S4symp} it follows that
\begin{align*}
\xi_2 \cdot y_2^TJ_nH^{-1}x_2 &=y_2^TJ_ny_2 = 0.
\end{align*}

Finally, we consider the three remaining entries in the first column,
\[
z^TJ_nH^{-1}y_2 = -(H^{-1}z)^TJ_ny_2
\]
 for $z = x_1, v_1, v_2.$ Due to $\vartheta_i H^{-1}v_i=u_i$ for $i =1,2,$ we obtain with \eqref{eq_S4symp}
\begin{align*}
\vartheta_i \cdot (H^{-1}v_i)^TJ_ny_2 = u_i^TJ_ny_2 = 0,
\end{align*}
while with \eqref{eq_y1} we have
\begin{align*}
\delta_1 \cdot (H^{-1}x_1)^TJ_ny_2 = y_1^TJ_ny_2 = 0.
\end{align*}
Hence, \eqref{eq_hinv4} holds.

\subsection{Step 5: $\operatorname{range}\{S_5\} = \mathcal{K}_{6}(H,u_1)+ \mathcal{K}_{4}(H^{-1},H^{-1}u_1)$  and Step 6:  $\operatorname{range}\{S_6\} = \mathcal{K}_{6}(H,u_1)+ \mathcal{K}_{6}(H^{-1},H^{-1}u_1)$}
We refrain from stating Steps 5 and 6 explicitly even so $u_2$ and $x_2$ are not displaying the general form of $u_k$ and $x_k$. This can only be seen from $u_3$ and $x_3$ which would be derived in Steps 5 and 6. As the derivations which lead to $u_3$ and $x_3$ are the same as in the general case for deriving $u_{k+1}$ and $x_{k+1},$ we directly proceed to the general case assuming that Algorithm \ref{alg} holds up to step $k.$

\subsection{Step 2k+1:  $\operatorname{range}\{S_{2k+1}\} = \mathcal{K}_{2k+2}(H,u_1)+ \mathcal{K}_{2k}(H^{-1},H^{-1}u_1)$}
Assume that we have constructed
\[
S_{2k} = [ y_k~ \cdots ~ y_1 ~ u_1 ~ \cdots ~ u_k \mid x_k~ \cdots ~ x_1 ~ v_1 ~ \cdots ~ v_k ] = [Y_k ~ U_k \mid X_k~V_k]\in \mathbb{R}^{2n \times 4k}
\]
such that $S_{2k}^TJnS_{2k} = J_{2k}$,
\begin{align}
H_{2k} &= J_{2k}^TS_{2k}^TJ_nHS_{2k} 
 = {\small
\begin{bmatrix}
-X_k^TJ_nHY_k & -X_k^TJ_nHU_k & -X_k^TJ_nHX_k & -X_k^TJ_nHV_k\\
-V_k^TJ_nHY_k & -V_k^TJ_nHU_k & -V_k^TJ_nHX_k & -V_k^TJ_nHV_k\\
Y_k^TJ_nHY_k & Y_k^TJ_nHU_k & Y_k^TJ_nHX_k & Y_k^TJ_nHV_k\\
U_k^TJ_nHY_k & U_k^TJ_nHU_k & U_k^TJ_nHX_k & U_k^TJ_nHV_k
\end{bmatrix}} \nonumber\\
& =
\begin{bmatrix}
0 & 0 & \Lambda_k & B_{kk}\\
0 & 0 & B_{kk}^T & T_k\\
\Delta_{k} & 0 & 0 & 0\\
0 & \Theta_{k} &0 & 0
\end{bmatrix}\label{eq_H2k}
\end{align}
 as in \eqref{eq_struct1} ($r=s=k$),
\begin{align}\label{eq_H2kinv}
 \tilde H_{2k} &= J_{2k}^TS_{2k}^TJ_nH^{-1}S_{2k}
 = {\small
\begin{bmatrix}
-X_k^TJ_nH^{-1}Y_k & -X_k^TJ_nH^{-1}U_k & -X_k^TJ_nH^{-1}X_k & -X_k^TJ_nH^{-1}V_k\\
-V_k^TJ_nH^{-1}Y_k & -V_k^TJ_nH^{-1}U_k & -V_k^TJ_nH^{-1}X_k & -V_k^TJ_nH^{-1}V_k\\
Y_k^TJ_nH^{-1}Y_k & Y_k^TJ_nH^{-1}U_k & Y_k^TJ_nH^{-1}X_k & Y_k^TJ_nH^{-1}V_k\\
U_k^TJ_nH^{-1}Y_k & U_k^TJ_nH^{-1}U_k & U_k^TJ_nH^{-1}X_k & U_k^TJ_nH^{-1}V_k
\end{bmatrix}} \nonumber\\
&=  \begin{bmatrix}
0 & 0 & \Delta_{k}^{-1} & 0\\
0 & 0 & 0 & \Theta_{k}^{-1}\\
E_{k} & G_{kk} & 0 & 0\\
G_{kk}^T&  F_{k}& 0&0
\end{bmatrix}
\end{align}
 as in \eqref{eq_myHinv}
and
\[
\operatorname{range}\{S_{2k}\} = \mathcal{K}_{{2k}}(H,u_1)+ \mathcal{K}_{{2k}}(H^{-1},H^{-1}u_1).
\]
The computational steps can be found in Algorithm \ref{alg}.

In this step the next two vectors $H^{2k}u_1$ and $H^{2k+1}u_1$ from $\mathcal{K}_{2k+2}(H,u_1)$ are added to the symplectic basis. Due to the previous construction, this is achieved by first  considering $Hv_k$.
$J$-orthogonalizing $ Hv_k$ against the columns of $S_{2k}$ yields
\begin{align*}
    w_u &=(I-S_{2k}J_2^TS_{2k}^TJ_n)Hv_k 
    = Hv_k -[X_k~~V_k~~-Y_k~~-U_k]\begin{bmatrix}
    Y_k^TJ_nHv_k\\
    U_k^TJ_nHv_k\\
    X_k^TJ_nHv_k\\
    V_k^TJ_nHv_k
    \end{bmatrix} \nonumber \\&
= Hv_k-\gamma_ky_k -\mu_ky_{k-1}-\beta_k u_{k-1} - \alpha_ku_k 
\end{align*}
as due to \eqref{eq_H2k}
\begin{align*}
 Y_k^TJ_nHv_k &=0,
  &  U_k^TJ_nHv_k=0,\\
    X_k^TJ_nHv_k&= \begin{bmatrix}
   -\gamma_k\\
-\mu_k\\
0\\
\vdots \\
0
    \end{bmatrix},
&    V_k^TJ_nHv_k = \begin{bmatrix}
0\\
\vdots\\
0\\
-\beta_k\\
    -\alpha_k
    \end{bmatrix}.
\end{align*}
Normalizing $w_u$ to length $1$  gives
\begin{equation}\label{eq_ukp1}
u_{k+1} = w_u/\chi_{k+1},
\end{equation}
where it is assumed that $\chi_{k+1} = \|w_u\|_2\neq 0.$

This step is finalized by $J$-orthogonalizing $Hu_{k+1}$ against the columns of $S_{2k}$
\begin{align*}
    w_v &=(I-S_{2k}J_{2k}^TS_{2k}^TJ_n)Hu_{k+1}
    = Hu_{k+1} -[X_k~~V_k~~-Y_k~~-U_k]\begin{bmatrix}
    Y_k^TJ_nHu_{k+1}\\
    U_k^TJ_nHu_{k+1}\\
    X_k^TJ_nHu_{k+1}\\
    V_k^TJ_nHu_{k+1}
    \end{bmatrix}.
\end{align*}
All entries of the last vector are zero.  The zeros in the first two blocks $Y_k^TJ_nHu_{k+1}$ and $U_k^TJ_nHu_{k+1}$   can be seen by using    $y_j = \delta_jH^{-1}x_j$ and $u_j=\vartheta_j H^{-1}v_j$ for $j = 1, \ldots, k$ as well as $H^{-T}J_nH=-J_n$:
\begin{align*}
   y_j^TJ_nHu_{k+1} /\delta_j &=  (H^{-1}x_j)^TJ_nHu_{k+1} =  -x_j^TJ_nu_{k+1} = 0,\\
  u_j^TJ_nHu_{k+1}/\vartheta_j  &= (H^{-1}v_j)^TJ_nHu_{k+1} =-v_j^TJ_nu_{k+1} = 0,
\end{align*}
due to the construction of $u_{k+1}$ as $J$-orthogonal against all columns of $S_{2k}.$

The zeros in the last block $V_k^TJ_nHv_k$ follow as $H$ is Hamiltonian and with
\[\chi_ju_j = Hv_j-\gamma_jy_j -\mu_jy_{j-1}-\beta_j u_{j-1} - \alpha_ju_j\]
for $j = 1, \ldots,k,$ (where we set $\beta_1=\mu_1 =0$ and $y_0=u_0=0$)
\begin{align}
  v_j^TJ_nHu_{k+1} &=v_j^T(J_nH)^Tu_{k+1} = -(Hv_j)^TJ_nu_{k+1}\nonumber\\
 &=-(\chi_ju_j+\gamma_jy_j+\mu_jy_{j-1}+\beta_ju_{j-1}+\alpha_ju_j)^TJ_nu_{k+1} = 0, \label{eq_hier}
\end{align}
again due to the construction of $u_{k+1}$ as $J$-orthogonal against all columns of $S_{2k}.$

With this we can show that the entries of the next to last block  $X_k^TJ_nHv_k$ are all zero.
First, with $\psi_1x_1 = H^{-1}u_1 -f_{11}v_1$ and $H^{-T}J_nH=-J_n$ we have
 \begin{align}
  \psi_1 \cdot x_1^TJ_nHu_{k+1} &=(H^{-1}u_1 -f_{11}v_1)^TJ_nHu_{k+1}
= u_1^T H^{-T}J_nHu_{k+1} -f_{11}v_1J_nHu_{k+1}\nonumber\\&
= -u_1^T J_nu_{k+1} -f_{11}v_1J_nHu_{k+1}
 =  0\label{eq_hier2}
\end{align}
due to the construction of $u_{k+1}$ as $J$-orthogonal against all columns of $S_{2k}$ and due to \eqref{eq_hier}.
Next, we use
\begin{align}\label{eq_xj}
\psi_{j+1} x_{j+1} =H^{-1}y_j -e_{jj}x_j-e_{j-1,j}x_{j-1}-g_{jj}v_{j}-g_{j1,j+1}v_{j+1}
\end{align}
 for $j = 1, \ldots, k-1$ (where we set $e_{01}=0$ and $x_0=0$, see Lines 16 and 29 of Algorithm \ref{alg})  for the other entries of the next to last block
 \begin{align*}
  \psi_{j+1} \cdot x_{j+1}^TJ_nHu_{k+1} &=-(e_{jj}x_j+e_{j-1,j}x_{j-1}+g_{jj}v_{j} +g_{j,j+1}v_{j+1})^TJ_nHu_{k+1} + y_j^TH^{-T}J_nHu_{k+1}\\
 &=-(e_{jj}x_j+e_{j-1,j}x_{j-1})^TJ_nHu_{k+1} -y_j^TJ_nu_{k+1}
\end{align*}
as $v_j^TJ_nHu_{k+1}=0$ due to \eqref{eq_hier}. Clearly, $y_j^TJ_nu_{k+1}=0$ by construction of $u_{k+1}.$ Thus, it remains to consider
 \begin{align*}
  \psi_{j+1} \cdot x_{j+1}^TJ_nHu_{k+1} &= -(e_{jj}x_j+e_{j-1,j}x_{j-1})^TJ_nHu_{k+1}.
\end{align*}
For $j=1$ we have with $x_0=0$ and \eqref{eq_hier2} that $ \psi_{2} \cdot x_{2}^TJ_nHu_{k+1} =0.$ With this, we get $\psi_{3} \cdot x_3^TJ_nHu_{k+1} =0,$ and, continuing in this fashion,
 \begin{align*}
  \psi_{j+1} \cdot x_{j+1}^TJ_nHu_{k+1} &= 0.
\end{align*}
Thus the expression for $w_v$ simplifies to
\[
w_v = Hu_{k+1}.
\]
Normalizing $w_v$ by $\vartheta_{k+1} = u_{k+1}^TJ_nHu_{k+1}$ to make sure it is $J$-orthogonal to $u_{k+1}$  yields
\begin{equation}\label{eq_vkp1}
v_{k+1} = Hu_{k+1}/\vartheta_{k+1}.
\end{equation}
Let
$S_{2k+1} = [y_k~\cdots~y_1~u_1~\cdots~u_{k+1}\mid x_k~\cdots~x_1~v_1~\cdots~v_{k+1}] 
= [Y_k~U_{k+1}\mid X_k~V_{k+1}]
\in \mathbb{R}^{2n \times 4k+2}.$
Then by construction
\begin{equation}\label{eq_Skp1symp}
S_{2k+1}^TJ_nS_{2k+1} = J_{2k+1}
\end{equation}
and
$\operatorname{range}\{S_{2k+1}\} = \mathcal{K}_{2k+2}(H,u_1)+ \mathcal{K}_{2k}(H^{-1},H^{-1}u_1).$
\subsubsection{The projected matrix $H_{2k+1} = J_{2k+1}^TS_{2k+1}^TJ_nHS_{2k+1}$}
Most of the entries in $H_{2k+1} = J_{2k+1}^TS_{2k+1}^TJ_nHS_{2k+1}$ (denoted in blue) are already known from $H_{2k}$ \eqref{eq_H2k}
\begin{align}
H_{2k+1} &=
{\tiny\left[\begin{array}{c|cc||c|cc}
{\color{blue}0} & {\color{blue}0} & -X_k^TJ_nHu_{k+1} & {\color{blue}\Lambda_k} & {\color{blue}B_{kk}} & -X_k^TJ_nHv_{k+1}\\ \hline
{\color{blue}0} & {\color{blue}0} &-V_{k}^TJ_nHu_{k+1} & {\color{blue}B_{kk}}^T & {\color{blue}T_k} & -V_k^TJ_nHv_{k+1}\\
-v_{k+1}^TJ_nHY_k & -v_{k+1}^TJ_nHU_k & -v_{k+1}^TJ_nHu_{k+1} & -v_{k+1}^TJ_nHX_k & -v_{k+1}^TJ_nHV_k & -v_{k+1}^TJ_nHv_{k+1}\\
\hline\hline
{\color{blue}\Delta_k} & {\color{blue}0} &Y_k^TJ_nHu_{k+1} & {\color{blue}0} & {\color{blue}0} & Y_k^TJ_nHv_{k+1}\\\hline
{\color{blue}0} & {\color{blue}\Theta_k} & U_{k}^TJ_nHu_{k+1} &{\color{blue}0} &{\color{blue}0} &U_k^TJ_nHv_{k+1}\\
u_{k+1}^TJ_nHY_k & u_{k+1}^TJ_nHU_k &u_{k+1}^TJ_nHu_{k+1} &u_{k+1}^TJ_nHX_k &u_{k+1}^TJ_nHV_k &u_{k+1}^TJ_nHv_{k+1}
\end{array}\right] }\nonumber\\
&=
\left[\begin{array}{c|cc||c|cc}
{\color{blue}0} & {\color{blue}0} & 0 & {\color{blue}\Lambda_k} & {\color{blue}B_{kk}} & \begin{smallmatrix}\mu_k\\0\\ \vdots\\0\end{smallmatrix}\\\hline
{\color{blue}0} & {\color{blue}0} &0 & {\color{blue}B_{kk}}^T & {\color{blue}T_k} & \begin{smallmatrix}0\\\vdots\\0\\\beta_{k+1}\end{smallmatrix}\\
0 & 0 & 0 & \begin{smallmatrix}\mu_k&0& \cdots&0\end{smallmatrix} & \begin{smallmatrix}0& \cdots&0&\beta_{k+1}\end{smallmatrix}& \alpha_{k+1}\\
\hline\hline
{\color{blue}\Delta_k} & {\color{blue}0} &0 & {\color{blue}0} & {\color{blue}0} & 0\\\hline
{\color{blue}0} & {\color{blue}\Theta_k} & 0 &{\color{blue}0} &{\color{blue}0} & 0\\
0 & 0 &\vartheta_{k+1} &0 &0 &0
\end{array}\right]
=
\left[\begin{array}{c|c||c|c}
0 & 0  & \Lambda_k & B_{k,k+1}\\ \hline
0 & 0  & B_{k,k+1}^T & T_{k+1}\\
\hline\hline
\Delta_k & 0  & 0 & 0\\\hline
0 & \Theta_{k+1}  &0 & 0
\end{array}\right]
.\label{eq_H2kp1}
\end{align}
The zeros in the third column (and hence in the last row) follow due to $Hu_{k+1} = \vartheta_{k+1} v_{k+1}$ and \eqref{eq_Skp1symp}.
The zeros in the first block $v_{k+1}^TJ_nHY_k$ of the third row follow due to $Hy_j=\delta_jx_j,$ for $j=1, \ldots, k,$ the ones in the second block $v_{k+1}^TJ_nHU_k$ due to   $Hu_j =\vartheta_j v_j,$  $j=1, \ldots,k.$ This also implies the zeros in the last row of the fourth and fifth block.

  Moreover, we obtain
\[
 v_{k+1}^TJ_nHV_{k-1} = 0
\]
from
\begin{align}\label{eq_uj}
\chi_{j+1}u_{j+1}= Hv_j-\gamma_jy_j-\mu_jy_{j-1}-\beta_ju_{j-1}-\alpha_ju_j
\end{align}
 for $j = 1, \ldots,k-1$ (where we set $\beta_0=\mu_0=0$ and $u_0=y_0=0,$ see Lines 11 and 24 in Algorithm \ref{alg}) as
\begin{align*}
v_{k+1}^TJ_nHv_{j} &=   v_{k+1}^TJ_n(\psi_{j+1}u_{j+1}+\gamma_jy_j+\mu_jy_{j-1}+\beta_ju_{j-1}+\alpha_ju_j)
=0
\end{align*}
due to the construction of $v_{k+1}$ as $J$-orthogonal against all columns of $S_{2k}.$
With this, $\psi_1x_1 =  H^{-1}u_1-f_{11}v_1$ and the recurrence for $x_{j}$ as in \eqref{eq_xj}, we observe that
\[
 v_{k+1}^TJ_nHX_{k-1} = 0
\]
holds. This can be seen step by step. Due to \eqref{eq_Skp1symp} and  $v_{k+1}^TJ_nHV_{k-1}=0,$ we have
\begin{align*}
\psi_1 v_{k+1}^TJ_nHx_{1} &= v_{k+1}^TJ_nH (H^{-1}u_1-f_{11}v_1)
=v_{k+1}^TJ_nu_1-    f_{11}v_{k+1}^TJ_nHv_1
=0,
\end{align*}
and with this and \eqref{eq_x2}  we have further
\begin{align*}
\psi_2 v_{k+1}^TJ_nHx_{2} &=-v_{k+1}^TJ_ny_{1}-e_{11}v_{k+1}^TJ_nHx_1 -g_{11}v_{k+1}^TJ_nHv_{1}-g_{12}v_{k+1}^TJ_nHv_{2}
=0.
\end{align*}
In this fashion we continue with the expression for $x_{j+1}$ as in Line 30 of Algorithm \ref{alg}  to obtain for $j = 2, \ldots, k-2$
\begin{align*}
\psi_{j+1} v_{k+1}^TJ_nHx_{j+1}&=-v_{k+1}^TJ_ny_{j}-e_{jj}v_{k+1}^TJ_nHx_j-e_{j-1,j}v_{k+1}^TJ_nHx_{j-1}\\&\qquad
-g_{jj}v_{k+1}^TJ_nHv_{j+1}-g_{j,j+1}v_{k+1}^TJ_nHv_{j}=0.
\end{align*}
Hence, \eqref{eq_H2kp1} holds.

\subsubsection{The projected matrix $ J_{2k+1}^TS_{2k+1}^TJ_nH^{-1}S_{2k+1}$}
Most of the entries in
\[\tilde H_{2k+1} = J_{2k+1}^TS_{2k+1}^TJ_nH^{-1}S_{2k+1}\]
 (denoted in blue) are already known from \eqref{eq_H2kinv}
\begin{align}
\tilde H_{2k+1} &=
{\tiny \left[\begin{array}{ccc|ccc}
{\color{blue}0} & {\color{blue}0 }& -X_k^TJ_nH^{-1}u_{k+1} &  {\color{blue}\Delta_k^{-1}} & {\color{blue}0}&-X_k^TJ_nH^{-1}v_{k+1}\\
{\color{blue}0} & {\color{blue}0} &-V_k^TJ_nH^{-1}u_{k+1} & {\color{blue}0} & {\color{blue}\Theta_k}^{-1} & -V_k^TJ_nH^{-1}v_{k+1}\\
-v_{k+1}^TJ_nH^{-1}Y_k & -v_{k+1}^TJ_nH^{-1}U_k & -v_{k+1}^TJ_nH^{-1}u_{k+1} & -v_{k+1}^TJ_nH^{-1}X_k & -v_{k+1}^TJ_nH^{-1}V_k & -v_{k+1}^TJ_nH^{-1}v_{k+1}\\
\hline
{\color{blue}E_k} & {\color{blue}G_k} &Y_k^TJ_nH^{-1}u_{k+1} &{\color{blue} 0} & {\color{blue}0} & Y_k^TJ_nH^{-1}v_{k+1}\\
{\color{blue}G_k}^T & {\color{blue}F_k} & U_k^TJ_nH^{-1}u_{k+1} &{\color{blue}0} &{\color{blue}0} & U_k^TJ_nH^{-1}v_{k+1}\\
u_{k+1}^TJ_nH^{-1}Y_k & u_{k+1}^TJ_nH^{-1}U_k &u_{k+1}^TJ_nH^{-1}u_{k+1} &u_{k+1}^TJ_nH^{-1}X_k &u_{k+1}^TJ_nH^{-1}V_k &u_{k+1}^TJ_nH^{-1}v_{k+1}
\end{array}\right] }\nonumber\\
&=
\left[\begin{array}{ccc|ccc}
{\color{blue}0} & {\color{blue}0} & 0 & {\color{blue}\Delta_k^{-1}} & {\color{blue}0} &0\\
{\color{blue}0} & {\color{blue}0} &0 & {\color{blue}0} & {\color{blue}\Theta_k}^{-1} &0\\
0 & 0 & 0 & 0 & 0 & 1/\vartheta_{k+1}\\
\hline
{\color{blue}E_k} & {\color{blue}G_k} &\begin{smallmatrix}g_{k,k+1}\\0\\\vdots\\0\end{smallmatrix} & {\color{blue}0} & {\color{blue}0}&0 \\
{\color{blue}G_k}^T & {\color{blue}F_{k}} & 0&{\color{blue}0} &{\color{blue}0}&0 \\
\begin{smallmatrix}g_{k,k+1}&0&\cdots&0\end{smallmatrix}  & 0 &f_{k+1,k+1} &0 &0 &0
\end{array}\right] \label{eq_hinvk}
=
\left[\begin{array}{cc|cc}
0 & 0  & \Delta_k^{-1} & 0 \\
 0 &0 & 0 & \Theta_{k+1}^{-1} \\
\hline
E_k & G_{k,k+1} & 0 & 0\\
G_{k,k+1}^T & F_{k+1} & 0&0
\end{array}\right]
\end{align}

With \eqref{eq_vkp1} and $H^{-T}J_nH=-J_n,$ we see that the entries $v_{k+1}^TJ_nH^{-1}z$ in the third block row are zeros (despite the last entry)
\begin{align}\label{eq_xx}
\vartheta_1 \cdot v_{k+1}^TJ_nH^{-1}z = u_{k+1}^TH^TJ_nH^{-1}z = -u_{k+1}^TJ_nz = 0
\end{align}
for $z \in\{y_1, \ldots, y_k, u_1, \ldots, u_{k+1},x_1,\ldots,x_k,v_1.\ldots,v_k\}$ due to the construction of $u_{k+1}$ as $J$-orthogonal to all columns of $S_{2k}.$ This implies the zeros in the last column of \eqref{eq_hinvk}. For the last entry we have
\begin{align*}
-\vartheta_{k+1} \cdot v_{k+1}^TJ_nH^{-1}v_{k+1} = -u_{k+1}^TH^TJ_nH^{-1}v_{k+1} =u_{k+1}^TJ_nv_{k+1} = 1.
\end{align*}
Thus, $-v_{k+1}^TJ_nH^{-1}v_{k+1}=1/\vartheta_{k+1}.$

With \eqref{eq_ukp1} the entries  in the upper part of the third column (as well as the entries in the fourth and fifth block of the last row) are zero as
\begin{align*}
\chi_{k+1}z^TJ_nH^{-1}u_{k+1} &=
 z^TJ_nH^{-1}(Hv_k-\gamma_ky_k-\mu_ky_{k-1}-\beta_ku_{k-1}-\alpha_ku_k)\\
&=z^TJ_nv_k - z^TJ_nH^{-1}(\gamma_ky_k+\mu_ky_{k-1}+\beta_ku_{k-1}+\alpha_ku_k)
=0
\end{align*}
for $z \in \{x_1, \ldots,x_k,v_1,\ldots,v_k\}$
due to \eqref{eq_xx} and \eqref{eq_H2kinv}.

The entries in $Y_{k-1}^TJ_nH^{-1}u_{k+1}$ are zero as $H^{-1}$ is Hamiltonian and \eqref{eq_xj} yield
\begin{align}\label{eq_yyy}
\begin{split}
y_j^TJ_nH^{-1}u_{k+1} &= -(H^{-1}y_j)^TJ_nu_{k+1} \\
&= -(\psi_{j+1}x_{j+1}-e_{jj}x_j-e_{j-1,j}x_{j-1}-g_{jj}v_j-g_{j,j+1}v_{j+1})^TJ_nu_{k+1}
=0
\end{split}
\end{align}
for $ j=1, \ldots, k-1$ due to the construction of $u_{k+1}$ as $J$-orthogonal to all columns of $S_{2k}.$

With this, we can show in a recursive manner that the entries in
\[
U_{k-1}^TJ_nH^{-1}u_{k+1} = -(H^{-1}U_{k-1})^TJ_nu_{k+1}
\]
 are zero by making use of $\chi_1 u_1 = \psi_1x_1-f_{11}v_1$ and \eqref{eq_uj}.
First we obtain
\begin{align}\label{eq_yy}
\chi_1 \cdot( H^{-1}u_1)^TJ_nu_{k+1} &=(\psi_1x_1-f_{11}v_1)^TJ_nu_{k+1}=0
\end{align}
due to the construction of $u_{k+1}$ as $J$-orthogonal to all columns of $S_{2k}.$
Next we observe
\begin{align*}
\chi_2 \cdot (H^{-1}u_2)^TJ_nu_{k+1}
&=(v_1-\gamma_1H^{-1}y_1-\alpha_1H^{-1}u_1)^TJ_nu_{k+1}=0,
\end{align*}
where the first term is zero as $u_{k+1}$ is $J$-orthogonal to $v_1,$ the second one due to \eqref{eq_yyy} and the third term due to \eqref{eq_yy}.  Continuing in this fashion, we have
\begin{align*}
\chi_j \cdot (H^{-1}u_j)^TJ_nu_{k+1}&=(v_{j-1}-\gamma_{j-1}H^{-1}y_{j-1}-\mu_{j-1}H^{-1}y_{j-2})^TJ_nu_{k+1}\\
&\qquad -(\beta_{j-2}H^{-1}u_{j-2}+\alpha_{j-1}H^{-1}u_{j-1})^TJ_nu_{k+1}
=0
\end{align*}
where the first term is zero as $u_{k+1}$ is $J$-orthogonal to $v_{j-1},$ the second and third one due to \eqref{eq_yyy}, and the fourth and fifth term due to the preceding observations.

Hence, \eqref{eq_hinvk} holds.

\subsection{Step 2k+2:  $\operatorname{range}\{S_{2k+2}\} = \mathcal{K}_{2k+2}(H,u_1)+ \mathcal{K}_{2k+2}(H^{-1},H^{-1}u_1)$}

Assume that we have constructed $S_{2k+1} =[Y_k~~U_{k+1}\mid X_k~~V_{k+1}] \in \mathbb{R}^{2n \times 4k+2}$ as in the previous section.

The two vectors $H^{-(2k+1)}u_1$ and $H^{-(2k+2)}u_1$ from $\mathcal{K}_{2k+2}(H^{-1},H^{-1}u_1)$ are added to the symplectic basis. Due to the previous construction, this is achieved by constructing $x_{k+1}$ from  $H^{-1}y_k$ and $y_{k+1}$ from $H^{-1}x_{k+1}.$
First  $ H^{-1}y_k$  is $J$-orthogonalized against the columns of $S_{2k+1}$:
\begin{align*}
    w_x &=(I-S_{2k+1}J_{2k+1}^TS_{2k+1}^TJ_n)H^{-1}y_k
    = H^{-1}y_k -[X_k~~V_{k+1}~~-Y_k~~-U_{k+1}]\begin{bmatrix}
    Y_k^TJ_nH^{-1}y_k\\
    U_{k+1}^TJ_nH^{-1}y_k\\
    X_k^TJ_nH^{-1}y_k\\
    V_{k+1}^TJ_nH^{-1}y_k
    \end{bmatrix}\nonumber\\
&= H^{-1}y_k -e_{kk}x_k-e_{k-1,k}x_{k-1}-g_{kk}v_{k}-g_{k,k+1}v_{k+1}
\end{align*}
due to \eqref{eq_hinvk}.
Normalizing $w_x$ to length $1$ gives
\begin{equation}\label{eq:xkp1}
x_{k+1}= w_x/\psi_{k+1},
\end{equation}
where we assume that $\psi_{k+1}=\|w_x\|_2 \neq 0.$

This step is finalized by $J$-orthogonalizing $H^{-1}x_{k+1}$ against the columns of $S_{2k+1}$:
\begin{align*}
    w_y&=(I-S_{2k+1}J_{2k+1}^TS_{2k+1}^TJ_n)H^{-1}x_{k+1}
    = H^{-1}x_{k+1} -[X_k~V_{k+1}\mid -Y_k~-U_{k+1}]\begin{bmatrix}
    Y_k^TJ_nH^{-1}x_{k+1}\\
    U_{k+1}^TJ_nH^{-1}x_{k+1}\\
    X_k^TJ_nH^{-1}x_{k+1}\\
    V_{k+1}^TJ_nH^{-1}x_{k+1}
    \end{bmatrix}.
\end{align*}
All entries $z^TJ_nH^{-1}x_{k+1} =-(H^{-1}z)^TJ_nx_{k+1}$ in the last vector are zero. In order to see this, let us first consider $z = v_j,  j = 1, \ldots, k+1.$ Due to $\vartheta_jv_j = Hu_j$, we have immediately
\[
(H^{-1}v_j)^TJ_nx_{k+1} = -u_j^TJ_nx_{k+1}/\vartheta_j =0.
\]
Next, we consider $z =x_j,  j = 1, \ldots, k,$ and make use of $\xi_{j}y_{j} = H^{-1}x_{j}$ to obtain
\[
(H^{-1}x_j)^TJ_nx_{k+1} = \zeta_jy_j^TJ_nx_{k+1} = 0.
\]
Rewriting \eqref{eq_xj} in terms of $H^{-1}y_j$, the case $z=y_j, j = 1, \ldots, k$ yields
\begin{align*}
(H^{-1}y_j)^TJ_nx_{k+1}& = (\psi_jx_j-e_{j-1,j-1}x_{j-1}-g_{j-1,j-1}v_{j-1}-g_{j-1,j}v_j)^TJ_nx_{k+1}
= 0.
\end{align*}
Finally, for $z =u_j$ we obtain from \eqref{eq_x1}
\[
(H^{-1}u_1)^TJ_nx_{k+1} = (\psi_1x_1-f_{11}v_1)^TJ_nx_{k+1}=0,
\]
from \eqref{eq_u2}
\[
(H^{-1}u_2)^TJ_nx_{k+1} = (v_1+\gamma_1H^{-1}y_1+\alpha_1H^{-1}u_1)^TJ_nx_{k+1}/\chi_2 =0,
\]
and from \eqref{eq_ukp1}
\begin{align*}
(H^{-1}u_{j+1})^TJ_nx_{k+1} &= (v_{j}+\gamma_{j}H^{-1}y_{j}+\mu_{j}y_{j-1}+\beta_{j}u_{j-1}+\alpha_jH^{-1}u_j)^TJ_nx_{k+1}/\chi_{j+1} =0
\end{align*}
for $j = 2, \ldots, k.$

Thus,
\begin{equation}\label{eq:ykp1}
y_{k+1} = H^{-1}x_{k+1}/\xi_{k+1}
\end{equation}
where we assume that
\[
\xi_{k+1} = (H^{-1}x_{k+1})^TJ_nx_{k+1} = x_{k+1}^TH^{-T}J_nx_{k+1}\neq 0.
\]
With the same argument as in \eqref{eq_xi} we see that
\[
    \delta _{k+1} = -\frac{1}{\xi_{k+1}}.
\]

Let $S_{2k+2} = [y_{k+1}~Y_k~u_{k+1}\mid x_{k+1}~X_k~V_{k+1}] \in \mathbb{R}^{2n \times 4k+4}.$
Then by construction $S_{2k+2}^TJ_nS_{2k+2} = J_{2k+2}$
and $\operatorname{range}\{S_{2k+2}\} = \mathcal{K}_{2k+2}(H,u_1)+ \mathcal{K}_{2k+2}(H^{-1},H^{-1}u_1).$
\subsubsection{The projected matrix $H_{2k+2} = J_{2k+2}^TS_{2k+2}^TJ_nHS_{2k+2}$}
Most of the entries in $H_{2k+2} = J_{2k+2}^TS_{2k+2}^TJ_nHS_{2k+2}$ (denoted in blue) are already known from \eqref{eq_H2kp1}
\begin{align}
&
{\footnotesize \left[\begin{array}{cc|c||cc|cc}
-x_{k+1}^TJ_nHy_{k+1} & -x_{k+1}^TJ_nHY_k & -x_{k+1}^TJ_nHU_{k+1}  & -x_{k+1}^TJ_nHx_{k+1} & -x_{k+1}^TJ_nHX_{k}& -x_{k+1}^TJ_nHV_{k+1} \\
-X_k^TJ_nHy_{k+1} & {\color{blue}0}  &{\color{blue}0} &-X_k^TJ_nHx_{k+1}  &  {\color{blue}\Lambda_k} & {\color{blue}B_{k,k+1}}\\ \hline
-V_{k+1}^TJ_nHy_{k+1} & {\color{blue}0} &{\color{blue}0} & -V_{k+1}^TJ_nHx_{k+1} &{\color{blue}B_{k,k+1}}^T & {\color{blue}T_{k+1}}\\
\hline\hline
y_{k+1}^TJ_nHy_{k+1} & y_{k+1}^TJ_nHY_k & y_{k+1}^TJ_nHU_{k+1}  & y_{k+1}^TJ_nHx_{k+1} & y_{k+1}^TJ_nHX_k & y_{k+1}^TJ_nHV_{k+1} \\
Y_k^TJ_nHy_{k+1}&{\color{blue}\Delta_k} &{\color{blue} 0} &Y_k^TJ_nHx_{k+1} &{\color{blue}0} & {\color{blue}0} \\\hline
U_{k+1}^TJ_nHy_{k+1}&{\color{blue}0} & {\color{blue}\Theta_{k+1}}  &U_{k+1}^TJ_nHx_{k+1}&{\color{blue}0} &{\color{blue}0}
\end{array}\right] }\nonumber\\
&=
\left[\begin{array}{cc|c||cc|c}
0 & 0 &  0 & \lambda_{k+1} & 0 & \begin{smallmatrix} 0 & \cdots& 0 & \gamma_{k+1}\end{smallmatrix}\\
0& {\color{blue}0}  & {\color{blue}0} & 0&{\color{blue} \Lambda_k} & {\color{blue}B_{k,k+1}}\\\hline
0&   {\color{blue}0} &{\color{blue}0} & \begin{smallmatrix} 0\\ \vdots \\0 \\ \gamma_{k+1}\end{smallmatrix}& {\color{blue}B_{k,k+1}}^T & {\color{blue}T_{k+1}}\\
\hline\hline
\delta_{k+1} & 0 & 0&0&0&0\\
0 &{\color{blue}\Delta_k} & {\color{blue}0} & 0 & {\color{blue}0} & {\color{blue}0}\\\hline
0 & {\color{blue}0} & {\color{blue}\Theta_{k+1}}  &0 &{\color{blue}0} & {\color{blue}0}
\end{array}\right]
=
\left[\begin{array}{cc||cc}
0  & 0&  \Lambda_{k+1} & B_{k+1,k+1}\\
0&   0 & B_{k1,k+1}^T & T_{k+1}\\
\hline\hline
\Delta_{k+1} & 0 & 0 & 0\\
 0 & \Theta_{k+1}  &0 &0
\end{array}\right]. \label{eq:H2kp2}
\end{align}
Making use of  \eqref{eq:ykp1} we obtain  $z^TJ_nHy_{k+1} = \delta_{k+1} z^TJ_nx_{k+1}=0$ for  all but two of the entries in the first column, that is, for
$z = x_1, \ldots, x_{k+1},$ $v_1,$ $\ldots,$ $v_{k+1},$ $y_1,$ $\ldots,$ $y_k,$ $u_1,$ $\ldots,$ $u_{k+1}.$ This gives the zeros in the fourth row as well.

For the entries $x_{k+1}^TJ_nHz$ in the first row we note that with \eqref{eq:xkp1} and $H^{-T}J_nH=J^T,$
\begin{align*}
x_{k+1}^TJ_nHz &= \psi_{k+1}(H^{-1}y_k-e_{kk}x_k-e_{k-1,k}x_{k-1}-g_{k-1,k}v_{k+1}-g_{kk}v_{k})^TJ_nHz
= 0
\end{align*}
for $z = y_1, \ldots, y_k, u_1, \ldots, u_{k+1}, x_1, \ldots, x_{k-2}, v_1, \ldots, v_{k-2}$ due to \eqref{eq_H2kp1} and  $S_{2k-1}^TJ_nS_{2k+1}=J_{2k+1}.$
Next, with $Hv_j = \chi_{j+1} u_{j+1} +\gamma_jy_j+\mu_jy_{j-1}+\beta_ku_{j-1}-\alpha_ju_j$ \eqref{eq_uj}, it follows for $j = k-1,k,k-1$ that
\begin{align*}
x_{k+1}^TJ_nHv_{k-2} =
x_{k+1}^TJ_nHv_{k-1} =
x_{k+1}^TJ_nHv_{k} =0
\end{align*}
as $S_{2k-1}^TJ_nS_{2k+1}=J_{2k+1}.$
With this and \eqref{eq:xkp1}  we obtain three more zero entries
\begin{align*}
\psi_{k-2}x_{k+1}^TJ_nHx_{k-2} =
\psi_{k-1}x_{k+1}^TJ_nHx_{k-1}=
\psi_{k}x_{k+1}^TJ_nHx_{k}
=0.
\end{align*}
This gives the zeros in the fourth column as well.

Hence, \eqref{eq:H2kp2} holds.

\subsubsection{The projected matrix $J^T_{2k+2}S_{2k+2}^TJ_nH^{-1}S_{2k+2}$}
Most of the entries in
\[\tilde H_{2k+2} = J^T_{2k+2}S_{2k+2}^TJ_nH^{-1}S_{2k+2}\]
(denoted in blue) are already known from \eqref{eq_hinvk},
\begin{align}
&
{\tiny \left[\begin{array}{ccc|ccc}
x_{k+1}^TJ_nH^{-1}y_{k+1} & {\color{red}0} & {\color{red}0}  & 1/\delta_{k+1} & {\color{red}0} & {\color{red}0}\\
X_k^TJ_nH^{-1}y_{k+1} & {\color{blue} 0} & {\color{blue}0} &{\color{red}0}  & {\color{blue}\Delta_k}^{-1} & {\color{blue}0} \\
V_{k+1}^TJ_nH^{-1}y_{k+1} & {\color{blue}0} & {\color{blue}0} & {\color{red}0}  &{\color{blue}0}& {\color{blue}\Theta_{k+1}}^{-1} \\
\hline
-y_{k+1}^TJ_nH^{-1}y_{k+1} &  -y_{k+1}^TJ_nH^{-1}Y_k&  -y_{k+1}^TJ_nH^{-1}U_{k+1} & -y_{k+1}^TJ_nH^{-1}x_{k+1} & -y_{k+1}^TJ_nH^{-1}X_k & -y_{k+1}^TJ_nH^{-1}V_{k+1}  \\
-Y_k^TJ_nH^{-1}y_{k+1}&{\color{blue}E_{k}} & {\color{blue}G_{k,k+1}}  &{\color{red}0} &{\color{blue}0 }& {\color{blue}0} \\
-U_{k+1}^TJ_nH^{-1}y_{k+1}&{\color{blue}G_{k,k+1}}^T & {\color{blue}F_{k+1}} &{\color{red}0}&{\color{blue}0} &{\color{blue}0}
\end{array}\right] }\nonumber\\
&=
\left[\begin{array}{cc|c||cc|c}
0 & {\color{red}0} & {\color{red}0}  & 1/\delta_{k+1} & {\color{red}0} & {\color{red}0} \\
0& {\color{blue}0} & {\color{blue}0}  & {\color{red}0}&  {\color{blue}\Delta_k}^{-1} & {\color{blue}0}\\\hline
0&  {\color{blue}0} & {\color{blue}0} & {\color{red}0}& {\color{blue}0} & {\color{blue}\Theta_{k+1}}^{-1} \\
\hline\hline
e_{k+1,k+1} &\begin{smallmatrix}e_{k,k+1} &0 & \cdots &0\end{smallmatrix} & \begin{smallmatrix}0&\cdots&0&g_{k+1,k+1}\end{smallmatrix} &0&0&0\\
\begin{smallmatrix}e_{k,k+1} \\ 0\\\vdots \\0\end{smallmatrix} &{\color{blue}E_{k}} & {\color{blue}G_{k,k+1}} & {\color{red}0} & {\color{blue}0} & {\color{blue}0}\\\hline
\begin{smallmatrix}0\\ \vdots\\0\\g_{k+1,k+1}\end{smallmatrix} & {\color{blue}G_{k,k+1}}^T& {\color{blue}F_{k+1}}  &{\color{red}0} &{\color{blue}0} & {\color{blue}0}
\end{array}\right]\label{eq:Hinv2kp2}\\
&=
\left[\begin{array}{cc||cc}
 0 & 0  &  \Delta_{k+1}^{-1} & 0\\
  0 & 0 & 0& \Theta_{k+1}^{-1} \\
\hline\hline
E_{k+1} & G_{k+1,k+1} & 0 & 0 \\
 G_{k+1,k+1}^T& F_{k+1}  &0 &0
\end{array}\right]\nonumber
\end{align}
All but one of the zeros in the first row (denoted in red) and the fourth column follow from the derivations in the previous section. Due to \eqref{eq:ykp1},
\[
x_{k+1}^TJ_nH^{-1}y_{k+1} = -(H^{-1}x_{k+1})^TJ_ny_{k+1} =-\xi_{k+1}y_{k+1}^TJ_ny_{k+1}=0,
\]
and the last zero in the first row/fourth columns follows.

Now, let us consider the first column. We have $X_k^TJ_nH^{-1}y_{k+1}=0$ as for $j=1, \ldots, k$
\[
x_j^TJ_nH^{-1}y_{k+1} = -(H^{-1}x_j)^TJ_ny_{k+1} = -\xi_jy_j^TJ_ny_{k+1}=0,
\]
and $V_{k+1}^TJ_nH^{-1}y_{k+1}=0$ as for $j=1, \ldots, k+1$
\[
v_j^TJ_nH^{-1}y_{k+1} = -(H^{-1}v_j)^TJ_ny_{k+1} = -u_j^TJ_ny_{k+1}/\vartheta_j=0.
\]
Next, observe that  $Y_{k-1}^TJ_nH^{-1}y_{k+1} =0$ as for $j = 1, \ldots, k-1$
\[
y_{j}^TJ_nH^{-1}y_{k+1} = \xi_j(H^{-1}x_j)^TJ_nH^{-1}y_{k+1}=\xi_jx_j^TJ_ny_{k+1} = 0.
\]
Finally, making use of \eqref{eq_uj} we observe that  $U_{k}^TJ_nH^{-1}y_{k+1} =0$ as for $j = 1, \ldots, k.$

Hence, \eqref{eq:Hinv2kp2} holds.



\end{document}